\documentclass[aop,MSNbibl,seceqn,citesort,dvips]{arximspdf}

\doi{10.1214/11-AOP730} %kopijuoti is PTS
\volume{41}
\issue{4}
\pubyear{2013}
\firstpage{2724}
\lastpage{2754}

\makeatletter

\setattribute{abstract}   {width}  {287pt}
\newcommand{\eqref}[1]{(\ref{#1})}

\def\lla{\langle}
\def\rra{\rangle}
\def\a{\alpha}

\def\cal{\mathcal}

\def\eps{\varepsilon}
\def\det{\operatorname{det}}

\def\Ee{\mathbb{E}}
\def\real{\mathbb{R}}
\def\definedas{\stackrel{\Delta}{=}}

\def\mbu{\mathbf{u}}
\def\bv{\mathbf{v}}
\def\calH{\mathcal{H}}
\def\calL{\mathcal{L}}
\def\calO{\mathcal{O}}

\newtheorem{theorem}{Theorem}[section]
\newtheorem{lemma}[theorem]{Lemma}
\newtheorem{proposition}[theorem]{Proposition}
\newproclaim{definition}[theorem]{Definition}
\newtheorem{cor}[theorem]{Corollary}
\newproclaim{remark}[theorem]{Remark}
\newproclaim{Example}[theorem]{Example}

\newcommand{\grad}{\nabla}
\newcommand{\minkmu}{\mathcal{M}^{\mu}}
\newcommand{\minkgk}{\mathcal{M}^{\gamma_k}}
\newcommand{\lips}{\mathcal{L}}
\newcommand{\argmin}{\operatorname{argmin}}
\newcommand{\support}{\mathcal{S}}
\newcommand{\Tube}{\operatorname{Tube}}
\def\overset{\stackrel}
\def\text{\mbox}

\makeatother

\begin{document}
\begin{frontmatter}

\title{Random fields and the geometry of Wiener space\thanksref{T1}}
\runtitle{Random fields and the Wiener space}
\thankstext{T1}{Supported in part by NSF Grants DMS-09-06801 and DMS-04-05970.}

\begin{aug}
\author[A]{\fnms{Jonathan E.} \snm{Taylor}\ead[label=e1]{jonathan.taylor@stanford.edu}}
\and
\author[B]{\fnms{Sreekar} \snm{Vadlamani}\corref{}\ead[label=e2]{sreekar@math.tifrbng.res.in}}
\runauthor{J. E. Taylor and S. Vadlamani}
\affiliation{Stanford University and TIFR-CAM}
\address[A]{Department of Statistics\\
Sequoia Hall\\
Stanford University\\
390 Serra Mall\\
Stanford, California 94305-4065\\
USA\\
\printead{e1}}

\address[B]{TIFR---Center for Applicable Mathematics\\
Sharadanagar, Chikkabommasandra\\
Post Bag 6503, GKVK Post Office\\
Bangalore 560 065\\
India\\
\printead{e2}}
\end{aug}

% HISTORY:
\received{\smonth{5} \syear{2011}}
\revised{\smonth{10} \syear{2011}}

% ABSTRACT
%
\begin{abstract}
In this work we consider infinite dimensional extensions of some finite
dimensional Gaussian geometric functionals called the
Gaussian Minkowski functionals.
These functionals appear as coefficients in the probability content of
a tube around a convex set $D \subset\real^k$ under
the standard Gaussian law $N(0, I_{k \times k})$.
Using these infinite dimensional extensions, we consider
geometric properties of some smooth random fields in the spirit of
[\textit{Random Fields and Geometry} (2007)
Springer] that can be expressed
in terms of reasonably smooth Wiener functionals.
\end{abstract}

% KEYWORDS
%
\begin{keyword}[class=AMS]
\kwd[Primary ]{60G60}
\kwd{60H05}
\kwd{60H07}
\kwd[; secondary ]{53C65}.
\end{keyword}

\begin{keyword}
\kwd{Wiener space}
\kwd{Malliavin calculus}
\kwd{random fields}.
\end{keyword}

\end{frontmatter}

%%-------------------------------------------------------------------------------------------------------
%%-------------------------------------------------------------------------------------------------------

%s1 #&#
\section{Introduction and motivation}
\label{secintro}

We start with a description of a certain
class of set functionals determined
by the canonical Gaussian measure on~$\real^k$.
By canonical,
we shall mean centered and having covariance $I_{k \times k}$.
Its density with respect to the Lebesgue measure on $\real^k$ is
therefore given
by $(2\pi)^{-k/2}e^{-\|x\|^2/2}$.
For this measure,
we consider computing the probability
content of a tube around $M$, leading us to a \textit{Gaussian tube
formula} which we state as
%
%e1 #&#
\begin{equation}
\label{eqfiniteGaussiantubeformula}
\gamma_k(M+\rho B_k) = \gamma_k(M)+\sum_{j=1}^{\infty}\frac{\rho
^j}{j!}\minkgk_j(M),
\end{equation}
where $\minkgk_j(M)$ is the $j$th \textit{Gaussian Minkowski Functional}
(GMF) of the set $M$.
If~$M$ is compact and convex, that is, if $M$ is a convex body, then we can
take the right-hand side \eqref{eqfiniteGaussiantubeformula} to be
a power
series expansion for the left-hand side. For certain~$M$, this expansion
must be taken to be a formal expansion, in the sense that up to terms of
some order, the left and right-hand side above agree. For example, if~$M$
is a centrally-symmetric cone such as the rejection region
for a $T$ or $F$ statistic, then $M$ has a singularity
at the origin in the sense that the geometric structure of the cone
around 0 is nonconvex and the expansion above is accurate only up to terms
of size~$O(\rho^{n-1})$.

Our interest in this tube formula lies in the appearance of these
coefficients in the expected Euler characteristic heu\-ristic \cite
{RFG,Worsley1,Worsley3,Worsley2}.

%%-----------------------------------------------------------------------------------------------------------------------------------------------------
%s1.1 #&#
\subsection{Expected Euler characteristic heuristic}

The Euler characteristic heu\-ristic was developed by Robert Adler and
Keith Worsley (cf., e.g.,~\cite{RFG,Worsley1,Worsley3,Worsley2})
to approximate the probability
\[
P\Bigl(\sup_{x \in M}f(x) \geq u\Bigr)
\]
with $E(\chi(A_t(f;M)))$, where $A_t(f;M)=\{x\in M\dvtx  f(x)\ge t\}\subset
M$, and $\chi$
is the Euler--Poincar\'e characteristic.

Let $M$ be an $m$-dimensional reasonably smooth manifold, with $(\xi
_1,\ldots,\xi_k)$ identically and
independently distributed copies of a Gaussian random field defined on
$M$. Subsequently, for any
$F\dvtx \real^k\to\real$, with two continuous derivatives, we can define
a new random field on $M$ given by
$f(x)=F(\xi_1(x),\ldots,\break\xi_k(x))$, for each $x\in M$.

Using the above Euler characteristic heuristic
for approximating the $P$-value for appropriately large values of $u$,
and Theorem $15.9.5$ of
\cite{RFG}, we have
%
%e2 #&#
\begin{eqnarray}
\label{eqP-value}
P\Bigl(\max_{x\in M}f(x)\ge u \Bigr)&\approx& E(\chi(A_u(f;M))
)
\nonumber
\\[-8pt]
\\[-8pt]
\nonumber
&=&\sum_{j=0}^m (2\pi)^{-j/2}\calL_j(M)\minkgk_j(F^{-1}[u,\infty)),
\end{eqnarray}
where $\minkgk_j(F^{-1}[u,\infty))$ for $j=0,1,\ldots$ are the GMFs
of the
set $F^{-1}[u,\infty)\subset\real^k$ that appear in (\ref{eqfiniteGaussiantubeformula}), and
$\calL_j(M)$ for $j=0,1,\ldots, m$ are the Lipschitz--Killing
curvatures (LKCs) of the manifold $M$
defined with respect to the Riemannian metric given by
%then the metric $g$ induced by a real valued Gaussian random field $
%$g(X,Y)=E(X\xi Y\xi),$ where $X\xi$ and $Y\xi$ are the directional
%derivatives of $\xi$.}
$g(X,Y)=E(X\xi_1 Y\xi_1)$, where $X$ and $Y$ are two vector fields on $M$,
with $X\xi_1$ and $Y\xi_1$ representing the directional derivatives
of $\xi_1$.

%%-----------------------------------------------------------------------------------------------------------------------------------------------------
%s1.2 #&#
\subsection{Curvature measures}

The LKCs for a large class
of subsets of any finite dimensional Euclidean space can be defined via
a \textit{Euclidean tube formula}.
In particular, let $M\subset\real^k$ be an $m$-dimensional set with
convex support cone, then writing
$\lambda_k$ as the standard $k$-dimensional Euclidean measure, $B_k$
as the $k$-dimensional
unit ball centered at origin, for small enough values of $\rho$, we have
%
%e3 #&#
\begin{eqnarray}
\label{eqLKCMink}
\lambda_k(M + \rho B_k)& =&
\sum_{j=0}^m \frac{\pi^{(n-j)/2}}{\Gamma({(n-j)}/{2}+1)}\rho
^{n-j}\calL_j(M)
\nonumber
\\[-8pt]
\\[-8pt]
\nonumber
& =&
\sum_{j=0}^m \frac{\rho^{n-j}}{(n-j)!} \theta_{n-j}(M),
\end{eqnarray}
where $\calL_j(M)$ is the $j$th LKC of the set $M$ with respect to the
usual Euclidean metric, and
$\theta_j(M)$'s are called the Minkowski functionals of the set $M$.

Geometrically, $\calL_{k-1}(M)$ for a smooth $(k-1)$-dimensional
manifold embedded in $\real^l$ is the $(k-1)$-dimensional Lebesgue
measure of the set $M$, and the other LKCs can be defined as
\[
\calL_{j}(M)=\frac{1}{s_{k-j}(k-1-j)!}
\int_{\partial M} P_{k-1-j}(\lambda_1(x),\ldots,\lambda_{k-1}(x))
\calH_{k-1}(dx),
\]
where $s_j$ is the surface area of a unit ball in $\real^j$, $(\lambda
_1(x),\ldots,\lambda_{k-1}(x))$ are the
principal curvatures at $x\in\partial M$,
and $P_i(\lambda_1(x),\ldots,\lambda_{k-1}(x))$ is the $i$th symmetric
polynomial in $(k-1)$ indices. In the case when the set $M$ is not unit
codimensional, then the definition
involves another integral over the normal bundle.

The (generalized) curvature measures defined this way
are therefore signed measures induced by the Lebesgue measure of the
ambient space. By
replacing the Lebesgue measure in \eqref{eqLKCMink} with an
appropriate Gaussian measure, we can
define a parallel Gaussian theory.
The GMFs in \eqref{eqfiniteGaussiantubeformula} play the
role of Minkowski functionals in the Gaussian theory.
In particular,
%
%e4 #&#
\begin{equation}
\label{eqminkfinite1}
\qquad\minkgk_j(M)\definedas(2\pi)^{-k/2}\sum_{m=0}^{j-1}
\pmatrix{ j-1 \vspace*{2pt}\cr m}
\Theta_{m+1}\bigl(M, H_{j-1-m}(\langle\eta,x\rangle)e^{-|x|^2/2}\bigr),
\end{equation}
{\spaceskip=0.18em plus 0.05em minus 0.02em where $\Theta_{m+1}(M,H_{j-1-m}(\langle\eta,x\rangle)e^{-|x|^2/2})$
is the integral of
$H_{j-1-m}(\langle\eta,\break x\rangle)e^{-|x|^2/2}$} with respect to the $(m+1)$th
generalized Minkowski curvature measure, and $H_k(y)$ is the $k$th
Hermite polynomial in $y$ (cf.~\cite{RFG}).

%%--------------------------------------------------------------------------------------------------------------------------------------------------------

%s1.3 #&#
\subsection{Our object of study: A richer class of random fields}

In this paper we intend to extend \eqref{eqP-value} to a larger class
of random fields $f$, which can be
expressed using $F\dvtx C_0[0,1]\to\real$, where $C_0[0,1]$ is the space
of continuous functions $f\dvtx [0,1]\to\real$,
such that $f(0)=0$, also referred to as the \textit{classical Wiener
space}, when equipped with the standard Wiener
measure on this sample space. In other words, we shall consider random
fields which can be expressed as some smooth \textit{Wiener functional}.
For instance, let us start with a smooth
manifold $M$ together with a Gaussian field $\{B^x(t)\dvtx t\in\real_+,
x\in M\}$ defined on it, such that its covariance
function is given by
%
%e5 #&#
\begin{equation}
\label{eqBcov}
E(B^x(t)B^y(s)) = s\wedge t C(x,y),
\end{equation}
where $C\dvtx M\times M\to\real$ is assumed to be a smooth
%of smoothness of $C$ will appear later in Section~\ref{secapp}, where
%we actually prove an extension of \eqref{eqP-value}}
function, with more details appearing in Section~\ref{secapp}, where
we actually prove an extension of \eqref{eqP-value}.
This infinite dimensional random field can be used to construct many more
random fields on~$M$, for instance, the following:
%
%ex1.1 #&#
\begin{Example}[(Stochastic integrals)]
\label{exstochint}
Let $V\dvtx \real\to\real$ be a smooth function, and consider the following
random field,
%
%e6 #&#
\begin{equation}
\label{eqstochintrf}
f(x)=\int_0^1 V(B^x(s)) \,dB^x(s) = F(B^x(\cdot)),
\end{equation}
where
$F\dvtx C_0 \rightarrow\real$ is the Wiener functional
\[
F(\omega) = \biggl(\int_0^1 V(B(t)) \,dB(t)\biggr)(\omega).
\]
This is clearly an extension of the random fields in \eqref{eqP-value}.
As a consequence of our extension of the Gaussian Minkowski
functionals to smooth Wiener functionals, we prove that, under suitable
smoothness conditions on $V$,
\[
\Ee(\chi(A_u(f;M)) ) = \sum_{j=0}^{\operatorname{dim}(M)}
(2\pi)^{-j/2} \lips_j(M) {\cal M}^{\mu}_j(F^{-1}[u,+\infty)).
\]
Our smoothness conditions
are rather strong in this paper: we assume
$V$ is $C^4$ with essentially polynomial growth. We need such strict
assumptions to ensure regularity of various conditional densities
derived from the random field \eqref{eqstochintrf} and its first
two derivatives at a point $x \in M$.
\end{Example}

A quick look at \eqref{eqP-value} reveals that in order to extend it
to the case when $F\dvtx C_0[0,1]\to\real$, we must
be able to define GMFs for infinite dimensional subsets of $C_0[0,1]$,
as $F^{-1}[u,\infty)\subset C_0[0,1]$.
In the present form, that is, \eqref{eqminkfinite1}, the definition
of GMFs appears to depend on
the summability of the principal curvatures of the set $\partial
(F^{-1}[u,\infty))$ at each point
$x\in\partial(F^{-1}[u,\infty))$ as well as the integrability of
these sums.
In infinite dimensions this summability requirement is equivalent
to an operator being trace class. This is quite a strong requirement,
and may be very
hard to check. Indeed, the natural summability requirements of operators
in the natural infinite dimensional calculus on $C_0$, the Malliavin
calculus, is the Hilbert--Schmidt class rather than the trace class.

Therefore, we shall first modify the definition of GMFs, from \eqref
{eqminkfinite1} to one which is more amenable
for an extension to the infinite dimensional case. This will be done in
Section~\ref{secprelimI}.

After setting up the notation and some technical background on the
Wiener space in Section
\ref{secprelimII}, the all important step, that of extending the
appropriate definition
of GMFs to the case of codimension one, \textit{smooth} subsets of the
Wiener space,
is accomplished in Section~\ref{sectubeformula}.
The characterization of GMFs in the infinite dimensional case
will be done precisely the same way as in the case of finite
dimensions, where, as noted earlier,
the GMFs are identified as the coefficients appearing in the Gaussian
tube formula.

Finally, in Section~\ref{secapp}, we use the infinite dimensional
extension of the GMFs to obtain an\vadjust{\goodbreak}
extension of \eqref{eqP-value}, for random fields which can be
expressed as stochastic integrals driven by
$B^x(\cdot)$ as defined in Example~\ref{exstochint}, and discuss
other possible implications of
the extension. Most of our methods in Section~\ref{secapp} are
invariant to the formulation of the random field
as a stochastic integral. Hence, should a random field satisfy all the
regularity conditions appearing in Section
\ref{secapp}, we expect our methods to work, modulo a few changes.

%%-------------------------------------------------------------------------------------------------------
%%-------------------------------------------------------------------------------------------------------
%s2 #&#
\section{Preliminaries I: The finite dimensional theory}
\label{secprelimI}

In this section we shall use the standard finite dimensional theory of
transformation of measure
for Gaussian spaces to modify the definition \eqref{eqminkfinite1}
of the GMFs to one which is more
suited to extension to the infinite dimensional case.

We begin by recalling some well-known facts about analysis on finite
\mbox{dimensional} Gaussian spaces
from Section $6.6.3$ of Chapter II of~\cite{Malliavin}, and Chapter~$3$ of~\cite{UstZakai00}.
Let $\gamma_k$ be the Gaussian measure on $\real^k$ given by\break $(2\pi
)^{-k/2}e^{-\|x\|^2/2}\,dx$, and $T$ a mapping from $\real^k$ into
itself, given by $T(x)=x+u(x)$, where $u\dvtx \real^k\to\real^k$ is
Sobolev differentiable
and $|u(x)-u(y)|\le c(\rho)|x-y|$ for any $x,y\in\real^k$ with
$|x-y|<\rho$. Then, the
Radon--Nikodym derivative of $\gamma_k\circ T$
with respect to the measure $\gamma_k$ is given by
%
%e7 #&#
\begin{equation}
\label{eqramerfinite}
\frac{d\gamma_k\circ T}{d\gamma_k}=|\det_2(I_{\real^k}+\nabla u)|
\exp\biggl(-\delta(u)-\frac{1}{2}\|u\|^2\biggr),
\end{equation}
where $\|\cdot\|$ is the usual Euclidean norm, and $\det_2$ is the
generalized Carleman--Fredholm
determinant.

Subsequently, for a smooth, unit codimensional, convex set $A\subset
\real^k$, let us define the tube
$\operatorname{Tube}(A,\rho)$ of width $\rho$ around the set $A$ as the set
$(A \oplus B(0,\rho))$, where $B(0,\rho)$
is the $k$-dimensional ball of radius $\rho$ centered at the origin.
Next, we shall define a signed distance function given by
\[
d_{\partial A}(x)= \cases{
\displaystyle\inf_{y\in\partial A}\|y-x\|, &\quad $\mbox{for } x\notin A,$\vspace*{2pt}\cr
\displaystyle-\inf_{y\in\partial A}\|y-x\|, & \quad$\mbox{for }
x\in\operatorname{Int}(A),$}
\]
where $\operatorname{Int}(A)$ denotes the interior of the set $A$.

Applying the co-area formula, and using the fact that $\|\nabla
d_{\partial A}\| =1$, we get
%
%e8 #&#
\begin{equation}
\label{eqminkfinite1b}
\gamma_k(\operatorname{Tube}(A,\rho)) = \gamma_k(A)+\int_0^{\rho}\int
_{d_A^{-1}(r)}
\frac{\exp(-\|x\|^2/2)}{(2\pi)^{n/2}} \,dx\, dr.
\end{equation}
For $r < \rho$ fixed, we can now use equation \eqref{eqramerfinite}
with any suitable transformation
$T_r\dvtx \real^k \rightarrow\real^k$ that agrees with
$x \mapsto x + r   \eta_x$
on $\{y\dvtx  d_A(y) \in(-\nu,\rho)\}$ for some small positive $\nu$.
Any such
transformation maps $\Tube(d_A^{-1}(r),\eps)$ to $\Tube(\partial A,
\eps)$
for $r < \rho$ and any $\eps< \nu$. Two further applications of the
co-area formula yield
\begin{eqnarray*}
&&\int_{d_A^{-1}(r)} \frac{\exp(-\|x\|^2/2)}{(2\pi)^{n/2}}
  \,dx\\
&&\qquad= \int_{\partial A} |\det_2(I_{\real^k}+r\nabla^2\,d_{\partial A})|
\exp\biggl(-r\delta(\nabla d_{\partial A})-\frac{1}{2}r^2\biggr) \frac
{\exp(-\|x\|^2/2)}{(2\pi)^{n/2}}   \,dx.
\end{eqnarray*}
Therefore, equation \eqref{eqminkfinite1b}
simplifies to
%
%e9 #&#
\begin{eqnarray}
\label{eqminkfinite2a}
&&\gamma_k(\operatorname{Tube}(A,\rho))\nonumber\\
&&\qquad = \gamma_k(A)+ \int_0^{\rho}\int
_{\partial A}
|\det_2(I_{\real^k}+r\nabla^2\,d_{\partial A})|\exp\biggl(-r\delta
(\nabla d_{\partial A})-\frac{1}{2}r^2\biggr)
\\
& &\hspace*{72pt}\qquad\quad{}  \times\frac{\exp(-\|x\|^2/2)}{(2\pi)^{n/2}} \,dx\, dr.\nonumber
\end{eqnarray}

Using a yet-to-be justified Taylor series expansion of the integrand
appearing in the above integral, we can
finally rewrite the GMFs as
%
%e10 #&#
\begin{equation}
\label{eqminkfinite2}
\qquad\minkgk_{j+1}(A) = \int_{\partial A}\frac{d^j}{d\rho^j}
\bigl(\det_2(I+\rho\nabla\eta)\exp\bigl(-\rho \delta(\eta)-\rho
^2/2\bigr)\bigr)\bigg|_{\rho=0}  \,da^{\partial A}(x),
\end{equation}
where $\eta=\nabla d_{\partial A}$ is the outward unit normal vector
field to the set $\partial A$, and $da^{\partial A}$ is the surface
measure of the set $\partial A$. Note that in the above expression we
have removed the modulus around the $\det_2$ part,
which can be justified by taking reasonably small values of $\rho$.
This new definition of GMFs involves terms which have obvious
extensions in the infinite dimensional case.

%%-------------------------------------------------------------------------------------------------------
%%-------------------------------------------------------------------------------------------------------
%s3 #&#
\section{Preliminaries II: The infinite dimensional theory}
\label{secprelimII}

In this section we recall some established concepts in
Malliavin calculus which we shall need in later sections.
We begin with an abstract Wiener space $(X,H,\mu)$, where $H$,
equipped with the inner product
$\langle\cdot,\cdot\rangle_H$, is a separable Hilbert space,
called the Cameron--Martin space, $X$ is a Banach space into which $H$
is injected continuously
and densely, and, finally, $\mu$ is the standard cylindrical Gaussian
measure on $H$. For the sake
of simplicity, one can appeal to the classical case when we have $H$ as
the space
of real-valued, absolutely continuous functions on $[0,1]$ with
$L^2([0,1])$ derivatives, which is
continuously embedded in $X=C_0([0,1])$ the space of real-valued
continuous functions $f$ on $[0,1]$,
such that $f(0)=0$.\vspace*{6pt}

%%-------------------------------------------------------------------------------------------------------
\textit{Sobolev spaces on Wiener space}

 Following the notation used in~\cite{Malliavin,Nualart-book,Watanabe93}, Sobolev spaces $D^p_{\a}(X;E)$ for $p>1$ and $\a>0$
are defined as the class of $E$ -valued functions $f\in L^p(X;E)$ such that
\[
\|f\|_{p,\a} \definedas\|(I-L)^{\a/2}f\|_{L^p(X;E)} <\infty,
\]
where $L$ is the Ornstein--Uhlenbeck operator defined on the Wiener
space. Writing~$D$ as the Gross--Sobolev derivative and $\delta$ as its dual under the
Wiener measure,
$L=-\delta D$.
The Sobolev spaces $D^p_{-\a}(X;\real)$ for $\a>0$ are the spaces of
distributions, defined as the dual of $D^q_{\a}(X;\real)$, where, as
usual, $p^{-1}+q^{-1}=1$.
Throughout this paper, whenever appropriate, we will adopt this convention.

The space of infinitely integrable, $\a$-smooth Wiener functionals is
given by
\[
D^{\infty-}_{\a}(X;\real)\definedas\bigcap_{1<p<\infty} D^p_{\a
}(X;\real).
\]
Consequently, let us define the analogous \textit{infinitely integrable}
random variables as
$L^{\infty-}(X;\real)\definedas D^{\infty-}_0(X;\real).$
Finally, we shall end this section with another definition which
translates to the regularity of
Wiener functionals.
%
%de3.1 #&#
\begin{definition}
\label{defnondegenerate}
For an $\real^k$-valued Wiener functional $F=(F_1,\ldots,F_k)$, the Malliavin
covariance (matrix) $\sigma^F=(\langle DF_i,DF_j\rangle_H)_{ij}$, and
the functional $F$ itself, is
called nondegenerate in the sense of Malliavin if $(\det\sigma
^F)^{-1}\in L^{\infty-}$, whenever
$\det\sigma^F$ is well defined.
\end{definition}%\vspace*{6pt}

%%-----------------------------------------------------------------------------------------------
\textit{H-Convexity}

In order to characterize the class of subsets of the Wiener space for
which we shall define the GMFs, we shall recall the notion of $H$-convexity.

%de3.2 #&#
\begin{definition} An $H$-convex functional is defined as a measurable
functional
$F\dvtx X\to\dvtx \real\cup\{\infty\}$ such that for any $h,k\in H$, $\a\in[0,1]$
%
%e11 #&#
\begin{equation}
\label{eqnH-convex}
F\bigl(\omega+\a h+(1-\a)k\bigr)\le\a F(\omega+h)+(1-\a)F(\omega+k)
\qquad\mbox{a.s.}
\end{equation}
\end{definition}

One of the properties of $H$-convex functionals which will be used in
later sections is that a necessary and sufficient condition for
a Wiener functional $F\in L^p$ for some $p>1$ to be $H$-convex is that
the corresponding $D^2 F$ must be a positive
and symmetric Hilbert--Schmidt operator valued distribution on $X$ (cf.
\cite{UstZakai00}).

%%-----------------------------------------------------------------------------------------------
%s3.1 #&#
\subsection{Quasi-sure analysis}
\label{subsecquasi}

In this section, most of which is based upon~\cite{Malliavin,Takeda},
we shall resolve some technical aspects of
defining integrals of Wiener functionals with respect to measures
concentrated on $\mu$-zero sets. Since all
Wiener functionals are \textit{de facto} defined up to $\mu$-zero sets,
thus, in order to be able to define the
integral of Wiener functionals with respect to measures which are
concentrated on $\mu$-zero sets,
we must resort to what is referred to as \textit{quasi-sure analysis},
which in turn relies on the concept of
\textit{capacities} on the Wiener space.

%de3.3 #&#
\begin{definition}
Let $1<p<\infty$ and $\a>0$. For an open set $O$ of $X$, we define
its $(p,\a)$-capacity
$C^p_{\a}(O)$ by
\[
C^p_{\a}(O)=\inf  \{\|U\|_{p,\a}\dvtx  U\in D^p_{\a}(X;\real),  U\ge
1  \ \mu\mbox{-a.e. on $O$}\}.
\]
For each subset of $A$ of $X$, we define its $(p,\a)$-capacity
$C^p_{\a}$ by
\[
C^p_{\a}(A)=\inf  \{C^p_{\a}(O)\dvtx  \mbox{$O$ is open and $O\supset
A$}\}.
\]
\end{definition}

These capacities are finer scales to estimate the size of sets in $X$
than~$\mu$. In particular,
a set of $(p,\a)$-capacity zero is always a $\mu$-zero set, but the
converse is not true in general.

A property $\pi$ is said to be true $(p,\a)$-quasi-everywhere (q.e.) if
\[
C^p_{\a}(\mbox{$\pi$ is not satisfied})=0.
\]

One of the most crucial steps in obtaining the \textit{co-area} formula
in the Wiener space,
which in turn is a necessary step to obtain the \textit{tube-formula} in
the Wiener space,
is to be able to extend ordinary Wiener functionals to sets of $\mu
$-zero measure.
Quasi-sure analysis lets us do precisely that and much more.

%de3.4 #&#
\begin{definition}
A measurable functional $F$ is said to have a $(p,\a)$-redefinition
$F^*$, satisfying
$F^* = F$ $\mu$-almost surely, and $F^*$ is $(p,\a)$-quasi-continuous,
if for all $\eps>0$, there exists an open set $O_{\eps}$ of $X$, such that
$C^p_{\a}(O_{\eps})<\eps$
and the restriction of $F^*$ to the complement set $O^c_{\eps}$ is
continuous under the norm
of uniform convergence on $X$.
\end{definition}

It can easily be seen that two redefinitions of the same functional
differ only on a set of
$(p,\a)$-capacity zero, thereby implying the uniqueness of a $(p,\a
)$-redefinition up to
$(p,\a)$-capacity zero sets. According to Theorem~$2.3.3$ of \cite
{Malliavin}, every functional
$F\in D^p_{\a}(X;\real^k)$ has a $(p,\a)$-quasi-continuous redefinition,
which can be taken to be in the first Baire class.

In what follows in the remainder of this section, we recall some
facts from the Malliavin calculus that will be helpful in our description
of a tube below. If $\a> 1$, one can make a statement similar to
Theorem $2.3.3$ of~\cite{Malliavin} related to the
differentiability of $F \in D^p_{\a}(X;\real)$, essentially
a form of Taylor's theorem with remainder.

%le3.5 #&#
\begin{lemma}
\label{lemmataylor}
Suppose $F \in D^p_{\a}(X;\real), \a> 1$. Then, for each $h \in H$
\[
\frac{1}{\varepsilon}\bigl(F(x+ \varepsilon h) - F(x) \bigr) -
\langle DF(x), h \rangle_H
\overset{D^{p_1}_{\a-1}(X;\real)}{\longrightarrow} 0
\]
for any $p_1 < p$.
\end{lemma}

\begin{pf}
Define
\begin{eqnarray*}
X_{n,h} &=& n\bigl(F(x+ h/n) - F(x) \bigr) \in D^{p_1}_{\a-1} \\
&=& \Lambda^{\a-1} Y_{n,h},\qquad   Y_{n,h} \in L^{p_1},\vadjust{\goodbreak}
\end{eqnarray*}
where $\Lambda= (I-L)^{-1/2}$ is the inverse of the Cauchy operator
\cite{Malliavin}. For each \mbox{$h \in H$}, $X_{n,h}$ converges in
$L^{p_1}$, so
the Kree--Meyer inequalities imply that $Y_{n,h}$ also converges in $L^{p_1}$.
A second application of the Kree--Meyer inequalities implies that
\[
\|Y_{n,h}-Y_{m,h}\|_{L^p} \approx\|X_{n,h}-X_{m,h}\|_{D^{p_1}_{\a-1}}.
\]
Or, $X_{n,h}$ is Cauchy in $D^{p_1}_{\a-1}$, hence, its limit $\langle
DF(x), h\rangle_H \in D^{p_1}_{\a-1}$.
\end{pf}

Hence, by the Borel--Cantelli property for the capacities $C^{p_1}_{\a
}$ (Corollary~IV.1.2.4 of
\cite{Malliavin}), for each $h \in H$ we can extract a sequence
$\varepsilon_n(h)$ such that
%
%e12 #&#
\begin{equation}
\label{eqdiffer}
\qquad C^{p_1}_{\a-1} \biggl(\biggl\{x\dvtx  \lim_{n \rightarrow\infty} \frac
{1}{\varepsilon_n(h)}
\bigl(F\bigl(x+\varepsilon_n(h)h\bigr)-F(x) \bigr) = \langle DF(x), h \rangle
\biggr\}^c \biggr) = 0.
\end{equation}

%co3.6 #&#
\begin{cor}
\label{corsmoothness}
Suppose $F \in D^p_{\a}(X;\real), \a> 1$ is nondegenerate and
$H_{\infty} \subset H$ is a countable dense subset.
Then,
%
%e13 #&#
\begin{eqnarray}
\label{eqsetsmoothness}
&&C^{p_1}_{\a-1} \biggl(\biggl\{x\dvtx  DF(x) \neq0\ \forall h
\in H_{\infty} \  \exists
\varepsilon_n(h) \rightarrow0 \mbox{ such that}
\nonumber
\\[-4pt]
\\[-8pt]
\nonumber
&&\hspace*{14pt}\qquad \lim_{n \rightarrow\infty}\biggl (\frac{1}{\varepsilon_n(h)}
\bigl(F^*\bigl(x+\varepsilon_n(h)h\bigr) - F^*(x) \bigr)
- \langle DF^*(x), h \rangle_H\biggr) = 0 \biggr\}^c \biggr) = 0.
\end{eqnarray}
\end{cor}

\begin{pf} The only thing that needs verifying beyond what was
pointed out above is that
$
C^{p_1}_{\a-1} ( \{x\dvtx  DF(x) \neq0 \}^c ) = 0.$
This follows from the Tchebycheff inequality (Theorem IV.2.2 of \cite
{Malliavin}) applied to
$\|DF\|_H \in D^{p_1}_{\a-1}$ and a Borel--Cantelli argument.
\end{pf}

%There is an obvious higher order version of Taylor's theorem above
%which we will use the second
%order version in our description of the tube below.
There is an obvious higher order Taylor expansion of
$F^{*}(x+\varepsilon_n(h)h)$, which we will use upto the second order term
in our description of the tube below.
If we are willing to sacrifice some moments, we can further specify in
Corollary~\ref{corsmoothness}
that the existence of the partial derivatives of $F$ as a limit at $x$
implies their existence as
limits at $x + h$ for all $h \in H_{\infty}, \|h\| \leq K$ for some
fixed, large $K$.

%co3.7 #&#
\begin{cor}
\label{coralltranslations}
Suppose $F \in D^p_{\a}(X;\real)$ and $H_{\infty} \subset H$ is a
countable dense subset.
Then, for all ${p_1}< p$
%
%e14 #&#
\begin{eqnarray}
\label{eqsetsmoothnessII}
&&
C^{p_1}_{\a-1} \biggl(\biggl\{x\dvtx  \forall h_1, h_2 \in
H_{\infty}, \|h_1\|
\leq K  \ \exists  \varepsilon_n(h_1, h_2) \rightarrow0 \mbox{ such that} \nonumber\\
&&\hspace*{15pt}\qquad  \lim_{n \rightarrow\infty} \biggl(\frac{1}{\varepsilon_n(h_1,h_2)}
\bigl(F^*\bigl(x+ h_1 + \varepsilon_n(h_1,h_2)h_2\bigr) - F^*(x + h_1)
\bigr)
\nonumber
\\[-8pt]
\\[-8pt]
\nonumber
& &\hspace*{176pt}{}\qquad  - \langle DF^*(x + h_1), h_2 \rangle_H \biggr) = 0
\biggr\} \biggr)\\
&&\qquad = 0.\nonumber
\end{eqnarray}
\end{cor}

\begin{pf} This follows from the fact that the translation operator
$
f(\cdot) \overset{T_h}{\mapsto} f(\cdot+h)
$
is a continuous map from $D^p_{\a}$ to $D^{p_1}_{\a}$ for any ${p_1}< p$
which follows directly from the Cameron--Martin theorem.
\end{pf}

%re3.8 #&#
\begin{remark}
\label{remlinearstructure}
Finally, we note that we can, by choosing $H_{\infty}$ appropriately,
choose the set, say,
$A$ in Corollary~\ref{coralltranslations}, in such a way that
$y \in A$ and $y + \varepsilon h \in A$ for all $h \in H_{\infty}$ and
for all $\varepsilon$ in some
countable dense subset of $\real$.
\end{remark}

%%-------------------------------------------------------------------------------------------------------
%%-------------------------------------------------------------------------------------------------------
%s4 #&#
\section{Key ingredients for a tube formula}
\label{sectube}

In this section we shall adopt a step-wise approach to reach our first
goal, that of
obtaining a (Gaussian) volume of the tube formula, for reasonably
smooth subsets of the Wiener space. The
three main steps are as follows: (i) characterizing subsets of the
Wiener space via Wiener functionals,
for which tubes, and thus GMFs, are well defined; (ii) assurance that
the surface measures
are well defined for the sets defined via the Wiener functionals; and
finally, (iii) a
change of measure formula for surface area measures corresponding to
the lower dimensional
surfaces of the Wiener space.

We shall first characterize the functionals for which the surfaces
measures are well defined,
subsequently, we shall prove a change of measure formula for the
surfaces defined via such
functionals. Finally, we shall define the class of sets for which the
tube formula and GMFs
are well defined by imposing more regularity conditions on the Wiener
functionals.

%%-------------------------------------------------------------------------------------------------------
%%-------------------------------------------------------------------------------------------------------
%s4.1 #&#
\subsection{The Wiener surface measures}
\label{subsecsurfacemeasure}

Let us start with a reasonably smooth, $\real^k$ valued Wiener
functional $F=(F_1,\ldots,F_k)$.
For $\mbu=(u_1,\ldots,u_k)\in\real^k$, we write $Z_{\mbu}=\bigcap
_{i=1}^kF_i^{-1}(u_i)$.
The sets $\{Z_{\mbu}\}_{{\mbu}\in\real^k}$ define a foliation of
hypersurfaces imbedded in $X$.

The surface measures of these foliations $Z_{\mbu}$ are closely
related to the density
$p_F$ of the push-forward measure $F_*(\mu)$ on $\real^k$ with
respect to the Lebesgue measure on $\real^k$.

Heuristically, writing $\delta_{\mbu}$ for the Dirac delta at ${\mbu
}\in\real^k$, the density $p_F$ can be defined as
%
%e15 #&#
\begin{equation}
\label{eqdeltadensity}
p_F({\mbu})=E(\delta_{\mbu}\circ F)
\end{equation}
as long as we can make sense of the composition $\delta_{\mbu}\circ F$.
For a smooth, real-valued Wiener functional $G$, we also expect the
following relation to hold:
%
%e16 #&#
\begin{equation}
\label{eqcondexp}
E[G \delta_{\mbu}\circ F]=E^{F={\mbu}}(G)\times p_F({\mbu}),
\end{equation}
where $E^{F={\mbu}}(G)$ is the conditional expectation of $G$ given
$F={\mbu}$,
assuming the composition $\delta_{\mbu} \circ F$ is well defined.

Making this heuristic calculation rigorous leads us back to the Sobolev
spaces of Section~\ref{secprelimII}, where
the object $(\delta_{\mbu} \circ F)$ is related to a \emph
{generalized Wiener functional}, that is, an element
of some $D^{p}_{-\a}$ for $p > 1, \a> 0$ through the pairing
\[
\langle G, \delta_{\mbu} \circ F \rangle_{D^q_{\a}, D^p_{-\a}} =
E[G \delta_{\mbu}\circ F],
\]
representing conditional expectation given $F=\mbu$ for any $G \in
D^{q}_{\a}$.
What is left to determine is, for a given $F$, which Sobolev spaces
contain $\delta_{\mbu} \circ F$.

The following theorem, the proof of which can be found in \cite
{Watanabe93}, provides the answer, taking us
one step closer to defining the surface measure corresponding to the
conditional expectation.

%th4.1 #&#
\begin{theorem}
\label{theoremwatanabe}
Let $F$ be an $\real^k$-valued, nondegenerate Wiener functional such that
$F\in D^{\infty-}_{1+\eps}(X;\real^k)$
for $\eps>0$, and the density $p_F$ of the law of $F$ is bounded.
Also, let
$0\le\beta<\min(\eps,\a)$ and $1<p<\infty$ satisfy
%
%e17 #&#
\begin{equation}
\label{eqalphabeta}
1<p<\frac{k}{\max\{(k+\beta-\min(\a,\eps)),0\}},
\end{equation}
and, finally, $\calO=\{z\in\real^k\dvtx  p_F(z)>0\}$. Then for $G\in
D^{q}_{\a}(X;\real)$,
with $\frac{1}{p}+\frac{1}{q}=1$, we have
%
%e18 #&#
\begin{equation}
\label{eqcondexpsobolev}
\zeta(u)= E(G \delta_u\circ F)\in W^{q}_{\beta}(\calO),
\end{equation}
where $W^{q}_{\beta}(\calO)$ is the Sobolev space of real-valued, weak
$\beta$-differentiable functions which are $q$-integrable.
\end{theorem}

Recall that for $F=(F_1,\ldots,F_k)\in D^{\infty-}_{1+\eps}$, the
density $p_F\in W^{\infty-}_{\eps}(\calO)$\vspace*{1pt}
(cf.~\cite{Nualart-book}). Now using the differentiability of the
density $p_F$ together with\vspace*{1pt}
equations \eqref{eqcondexp}, \eqref{eqcondexpsobolev} and the
algebraic structure of
the Sobolev spaces, we have $E^F(G)\in W^{q}_{\beta}(\calO),$
for any $G\in D^{q}_{\a}(X;\real)$, where $(E^F(G))(\mbu
) = E^{F=\mbu}(G)$.

That is, for each $F \in D^{\infty-}_{1+\eps}(X;\real)$, there
exists a continuous mapping
$E^F\dvtx$ $D^q_{\a}(X;\real) \rightarrow W^q_{\beta}(\calO).$ This, in
turn, induces a dual map\vadjust{\goodbreak}
$(E^F)^*\dvtx W^p_{-\beta}(\calO) \rightarrow D^p_{-\a}(X;\real)$
defined via the dual relationship
%
%e19 #&#
\begin{equation}
\label{eqwatanabemap}
\langle E^F(G),v\rangle_{W^q_{\beta}(\calO),W^p_{-\beta}(\calO)} =
\langle G,(E^F)^*v\rangle_{D^q_{\a}(X;\real),D^p_{-\a}(X;\real)}.
\end{equation}
Informally, this map, sometimes referred to as the Watanabe
map (see Section 6 of Chapter III of~\cite{Malliavin}), is just
composition, that is,
$(E^F)^*v = v \circ F.$

The object $(E^F)^*\delta_{\mbu}$ is almost the surface measure
needed in
\eqref{eqminkfinite2}, but it is just a generalized Wiener
functional, that is, distribution on $X$, at this point.
If we are to justify our Taylor series expansion via a dominated
convergence argument, we need to know that it
has a representation as a measure on $X$.

Clearly, for positive $G\in D^{q}_{\a}(X;\real)$, we shall have
\[
\langle G,(E^F)^*\delta_{\mbu}\rangle_{D^q_{\a}(X;\real),D^p_{-\a
}(X;\real)} = E^{F={\mbu}}(G)> 0.
\]
Therefore, $(E^F)^*\delta_{\mbu}\in D^p_{-\a}(X;\real)$ defines a
\textit{positive generalized Wiener functional}. Next, Theorem 4.3 of
\cite{Sugita} together with the conditions stated in
Theorem~\ref{theoremwatanabe}
implies that for each $u\in\calO$, there exists a finite positive
Borel measure $\nu^{F,{\mbu}}$
defined on Borel subsets of the Wiener space $X$, supported on
$F^{-1}(\mbu)$, such that
\[
E^{F={\mbu}}(G)=\int_X G^*(x) \nu^{F,{\mbu}}(dx)
\]
for all $G\in D^{q}_{\a}(X;\real)$, with $G^*$ its $(q,\a)$-quasi
continuous redefinition.

The measure $\nu^{F,{\mbu}}$ defined is a probability measure on the
set $F^{-1}({\mbu})$.
Using Airault and Malliavin's arguments in~\cite{AiraultMalliavin},
an appropriate area
measure $da^{Z_{\mbu}}$, corresponding to the measure $\nu^{F,{\mbu
}}$, can be defined as
%
%e20 #&#
\begin{equation}
\label{eqsurfacemeasure}
\int_X G^*(x) \,da^{Z_{\mbu}}(x)\definedas p_F({\mbu})\int
G^*(x)(\det(\sigma_F))^{1/2} \nu^{F,{\mbu}}(dx),
\end{equation}
where $\sigma_F$ is the Malliavin covariance matrix. Note that
the surface measure depends
only on the geometry of the set $Z_{\mbu}$, whereas the conditional
probability measure depends on the functional
from which the set is derived, thus the superscripts on the respective measures.
We are now in a position to justify at least part of \eqref{eqminkfinite2}.

%th4.2 #&#
\begin{theorem}
\label{theoremexistence}
Let $F$ be a $\real$-valued nondegenerate Wiener functional
such that $F \in D^{\infty-}_{2+\eps}(X;\real^k)$
and the density $p_F$ of the law of $F$ is bounded. Define the
unit normal vector field $\eta= DF/\|DF\|_H$.
Furthermore, suppose that:
\begin{itemize}
\item$E(\exp(\rho  \delta(\eta)) ) <
\infty$ for $\rho$ in some neighborhood of 0;

\item$E(\exp(\rho^2 \|D\eta\|^2_{\otimes^2 H} )
) < \infty$ for $\rho$ in some
neighborhood of 0.
\end{itemize}
Then, for $0 \leq\rho< \rho_c$ for some nonzero critical radius
%
%e21 #&#
\begin{eqnarray}
\label{eqtaylorseries}
\hspace*{8pt}&&
\int_0^{\rho}\int_{F^{-1}(\mbu)}
\det_2(I_H+r D\eta)\exp\biggl(-r\delta(\eta)-\frac
{1}{2}r^2\biggr) \,da^{Z_{\mbu}}  \,dr
\nonumber
\\[-4pt]
\\[-12pt]
\nonumber
&&\qquad= \sum_{j \geq1} \frac{\rho^j}{j!} \int_{F^{-1}(\mbu)}\frac
{d^{j-1}}{dr^{j-1}}
\bigl(\det_2(I_H+r D\eta)\exp\bigl(-r   \delta(\eta)-r^2/2\bigr)\bigr)
\bigg|_{r=0} \,da^{Z_{\mbu}}.
\end{eqnarray}
\end{theorem}

Before proving the above theorem, we shall state a few results
concerning the regularity of
functions of smooth Wiener functionals.

%pr4.3 #&#
\begin{proposition}
\label{propregularity}
Let $\a>0$ and $U\in D^{\infty-}_{\a}(X;\real)$.
\begin{itemize}
\item If $\exp(U)\in L^p(X;\real)$, then $\exp(U)\in D^{p'}_{\a'}$
where $p'=p^2(\a-\a')$ for $\a'<\a$.
\item If $U>0$ $\mu$ almost surely and $1/U\in L^p$, then $1/U\in
D^{p'}_{\a'}$ where $p'=p^2(\a-\a')$ for $\a'<\a$.
\end{itemize}
\end{proposition}

We shall skip the proofs of the above, as these can be proved by
replicating the proofs of Theorems $1.4$ and $1.5$ of Watanabe \cite
{Watanabe93}.

\begin{pf*}{Proof of Theorem \protect\ref{theoremexistence}} This is just
dominated convergence combined with the
nondegeneracy of $F$ as well as the following bound (cf. Theorem~9.2 of
\cite{Simon}):
\[
|\det_2(I+A)| \leq\exp(C \|A\|_{\otimes^2 H}^2 )
\]
for some fixed $C > 0$.

Note that while using the dominated convergence, we are inherently
assuming the \textit{well definedness} of integrals
of $\exp(\rho  \delta(\eta))$ and $\exp(\rho
^2 \|D\eta\|^2_{\otimes^2 H})$ with respect to the surface
measure $da^{Z_{\mbu}}$, which requires
%
%e22 #&#
\begin{eqnarray}
\label{eqexpregularity1}
\exp(\rho  \delta(\eta))
&\in& D^{q}_{\a}(X;\real)\qquad  \mbox{such that }q> 1/(\min\{\a,1+\eps
\}),
\\
\label{eqexpregularity2}
\quad\exp(\rho^2 \|D\eta\|^2_{\otimes^2 H})
&\in& D^{q}_{\a}(X;\real)  \qquad\mbox{such that }q> 1/(\min\{\a,1+\eps
\}).
\end{eqnarray}
Now, using Theorem $1.5$ of~\cite{Watanabe93}, we have
$\eta\in D^{\infty-}_{1+\eps'}$ for all $\eps'<\eps$.
Subsequently, using the above proposition together with the
assumption involving the existence of exponential moments, we have
$\exp(\rho  \delta(\eta)), \exp(\rho^2 \|
D\eta\|^2_{\otimes^2 H})\in D^{p}_{\eps''}(X;\real),$
such that $p=(\frac{\rho_c}{\rho})^2(\eps'-\eps''),$ where $\eps
''<\eps'$. In order to satisfy \eqref{eqexpregularity1} and
\eqref{eqexpregularity2}, we must choose $\eps'$ and $\eps''$
such that $\rho<\rho^2_c\eps''(\eps'-\eps'')$.
\end{pf*}

%re4.4 #&#
\begin{remark}
Note that Theorem~\ref{theoremexistence} does not say that the
Gaussian measure of the tube is given by the power series in \eqref
{eqtaylorseries}. Rather, it gives
conditions on the sets $Z_{\mbu}=F^{-1}(\mbu)$
for which the coefficients in the power series are well defined.
These conditions allow us to define GMFs for level sets of functions
that are not necessarily $H$-convex. However, for such functions we will
lose the interpretation of the power series in \eqref
{eqtaylorseries} as
an expansion for the Gaussian measure of the tube. This is similar
to the distinction between the formal and exact versions of the Weyl/Steiner tube formulae
\cite{TakemuraKurikiEquivalence2002}.
\end{remark}

%%----------------------------------------------------------------------------------------------------

%s4.2 #&#
\subsection{Change of measure formula: A Ramer type formula for
surface measures}
\label{subsecramer}

After assuring ourselves of the existence of the surface Wiener
measures, we shall now move
onto proving a change of measure formula for the surface measures given
by equation
\eqref{eqsurfacemeasure}.

To begin with, let $F\in D^{\infty-}_{1+\eps}(X;\real^k)$
so that we can define the surface measure using Theorem~\ref{theoremwatanabe}.
In order to obtain a change of measure formula for the
lower-dimensional subspaces
of the Wiener space, we shall start with the standard change of measure formula
on the Wiener space $X$.
Let us define a mapping $T_{\eta}\dvtx X\to X$ given by $T_{\eta
}(x)=x+\eta_x$,
for some smooth $\eta\dvtx X\to H$. Moreover, let $U$ be an open subset of
$X$, and:
\begin{longlist}[(1)]
\item[(1)]$T_{\eta}$ is a homeomorphism of $U$ onto an open subset of $X$,
\item[(2)]$\eta$ is an $H$-valued $C^1$ map and its $H$ derivative at each
$x\in U$
is a Hilbert--Schimdt operator on $H$.
\end{longlist}
This transformation induces two types of changes on the initial measure
$\mu$ defined on $X$.
These two induced measures can be expressed as
\begin{eqnarray*}
P(A) & = & \mu(T^{-1}_{\eta}(A)) = T_{\eta}^*\mu(A),\\
Q(A) & = & \mu(T_{\eta}(A)) = (T_{\eta}^{-1})^*\mu(A)
\end{eqnarray*}
for $A$ a Borel set of $X$.

Ramer's formula for change of measure on $X$, induced by a
transformation defined on $X$ and
satisfying the above conditions, gives an expression for the
Radon--Nikodym derivative
of $\mu\circ T_{\eta}$ with respect to $\mu$ and can be stated as follows:
%
%e24 #&#
\begin{equation}
\label{eqRamer}
\frac{dQ}{d\mu} = \bigl|\det_2\bigl(I_H+\nabla\eta(x)\bigr)\bigr|\exp\biggl(-\delta(\eta)
-\frac{1}{2}\|\eta(x)\|^2_H\biggr)
\definedas    Y_{\eta}(x),
\end{equation}
where $\delta(\eta)$ denotes the Malliavin divergence of an
$H$-valued vector
field $\eta$ in $X$. The proof of this result can be found in \cite
{Ramer,UstZakai00}.
It is to be noted here that, for appropriately smooth transformations,
a similar result for
$d(\mu\circ T^{-1}_{\eta})/d\mu$ can be obtained by using the
relationship between
$d(\mu\circ T_{\eta})/d\mu$ and $d(\mu\circ T^{-1}_{\eta})/d\mu$
given by
\[
\frac{d\mu\circ T_{\eta}}{d\mu}(x)
=\biggl(\frac{d\mu\circ T^{-1}_{\eta}}{d\mu}(T_{\eta}x)\biggr)^{-1}.
\]

The following theorem is the first step toward obtaining similar
formulae for change of measure
on lower-dimensional subsets of the Wiener space.

%th4.5 #&#
\begin{theorem}
\label{theoremsurfacemeasure}
Let $F\in D^{\infty-}_{1+\eps}(X;\real^k)$ satisfy the conditions
from Theorem~\ref{theoremwatanabe}, and $\a,\beta$ and $p$ be as given in \eqref
{eqalphabeta}.
Then, there exists a sequence of probability measures $\{\nu^{F,\mbu
}_n\}_{n\ge1}$
defined on Borel subsets of $X$ such that the measures $\{\nu^{F,\mbu
}_n\}_{n\ge1}$ are absolutely continuous with
respect to the Wiener measure $\mu$ and
the sequence $\{\nu^{F,\mbu}_n\}_{n\ge1}$ converges weakly to $\nu
^{F,\mbu}$.
\end{theorem}

\begin{pf} Let us choose a sequence of positive distributions on
$\calO$ given by
$\{v_n\}_{n\ge1}\subset W^p_{-\beta}(\calO)$, such that it converges
to $\delta_{\mbu}$
weakly in $W^p_{-\beta}(\calO)$ and that $\int v_n(\xi) \,d\xi=1$
for all $n\ge1$.
Then define the measures $\nu^{F,\mbu}_n$ as
\[
\int G(x) \nu^{F,\mbu}_n(dx)=\frac{1}{p_F(\mbu)}\int G(x)
v_n(F(x))  \mu(dx)
\]
for all measurable $G$ on $(X,\mu)$.
In view of \eqref{eqwatanabemap},
we can clearly identify the restriction of the measures $\nu^{F,\mbu
}_n$ to
$D^q_{\a}(X;\real)$ with $(E^F)^*v_n$. Now, by construction, $\{
(E^F)^*v_n\}_{n\ge1}$
converges to $(E^F)^*\delta_{\mbu}$ in $D^p_{-\a}(X;\real)$ and
their limit is a
nonnegative generalized Wiener functional. Therefore, using Lemma 4.1
of~\cite{Sugita},
we see that the measures $\nu^{F,\mbu}_n$ converge weakly to $\nu
^{F,\mbu}$.
\end{pf}

Thus, the surface (probability) measure of $Z_{\mbu}$, or the conditional
probability measure corresponding to $\{F=\mbu\}$, for any $\mbu\in
\calO$
can also be defined as
%
%e25 #&#
\begin{eqnarray}
\label{eqsurfaceprobmeasure}
\int G^*(x) \nu^{F,\mbu}(dx)&=&\lim_n \int G^*(x) \nu^{F,\mbu
}_n(dx)
\nonumber
\\[-8pt]
\\[-8pt]
\nonumber
& =&
\lim_n \frac{1}{p_F(\mbu)}\int G^*(x)v_n(F(x)) \mu(dx)
\end{eqnarray}
for the appropriate class of Wiener functionals $G$, which, as noted
earlier, depends on the regularity of $F$.

Let us now define a mapping $T_{\rho,\eta}\dvtx X\to X$ given by $T_{\rho
,\eta}(x)=x+\rho\eta_x$,
for some $\eta\in D^{\infty-}_{1+\eps}(X;H)$. We shall study the
change that the mapping
$T_{\rho,\eta}$ induces on the surface measure of $F^{-1}(\mbu)$.
Note that
\[
T_{\rho,\eta}(F^{-1}(\mbu))=\{x+\rho\eta_x\dvtx  x\in F^{-1}(\mbu)\}.
\]
Set $F_{\rho,\eta}=F\circ T_{\rho,\eta}^{-1} \in D^{\infty
-}_{1+\eps}$, so that
\[
F_{\rho,\eta}^{-1}(\mbu)=T_{\rho,\eta}(F^{-1}(\mbu)) \definedas
Z_{\mbu}^{\eta,\rho}.
\]
Using the above theorem and the relationship \eqref
{eqsurfacemeasure}, the area measure for $Z_{\mbu}^{\eta,\rho}$
can now
be identified as
\[
\int_{Z_{\mbu}^{\eta,\rho}}G^*\,da^{Z_{\mbu}^{\eta,\rho}} = \lim
_n \int_{X}G^*(y)[\det\sigma_{F_{\rho,\eta}}(y)]^{1/2}
v_n(F_{\rho,\eta}(y)) \mu(dy).
\]

Now using the transformation $y=T_{\rho,\eta}(x)$ and replacing the
function $G^*(\cdot)$
by $G^*(T_{\rho,\eta}^{-1}(\cdot))$ and, finally, using the standard
Ramer's\vadjust{\goodbreak} formula from
equation \eqref{eqRamer}, we get
%
%e26 #&#
\begin{eqnarray}
\label{eq1}
&&
\int_{Z_{\mbu}^{\eta,\rho}}G^*(T^{-1}_{\rho,\eta
}y)\,da^{Z_{\mbu}^{\eta,\rho}}\nonumber\\
&&\qquad = \lim_n \int_{X}G^*(x)[\det\sigma_{F_{\rho,\eta}}(T_{\rho
,\eta}x)]^{1/2}
v_n(F_{\rho,\eta}(T_{\rho,\eta}x)) Y^{\eta}_{\rho}(x) \mu(dx) \\
&&\qquad = \lim_n \int_{X}G^*(x)[\det\sigma_{F_{\rho,\eta}}(T_{\rho
,\eta}x)]^{1/2}
v_n(F(x)) Y^{\eta}_{\rho}(x) \mu(dx),\nonumber
\end{eqnarray}
where $Y^{\rho,\eta}(x)$ is the Radon--Nikodym derivative of the
measure $\mu\circ T_{\rho,\eta}$
with respect to the measure $\mu$ and, as a result of \eqref
{eqRamer}, can be expressed as
%
%e27 #&#
\begin{equation}
\label{eqRamerrho}
Y^{\eta}_{\rho}(x)= \bigl|\det_2\bigl(I_H+\rho\nabla\eta(x)\bigr)\bigr| \exp\bigl(-\rho
\delta(\eta)
-\tfrac{1}{2}\rho^2\|\eta\|^2_H\bigr).
\end{equation}

Using the definitions of $F_{\rho}$, the surface (probability and
area) measures
and, finally, rearranging the terms, we can rewrite \eqref{eq1} as
%
%e28 #&#
\begin{eqnarray}
\quad&&
\int_{Z_{\mbu}^{\eta,\rho}}G^*(T^{-1}_{\rho,\eta
}y)\,da^{Z_{\mbu}^{\eta,\rho}}\nonumber\\
&&\qquad= \lim_n \int_{X}G^*(x)[\det\sigma_{F_{\rho,\eta}}(T_{\rho,\eta
}x)]^{1/2}
Y^{\eta}_{\rho}(x) v_n(F(x)) \mu(dx)\nonumber\\
&&\qquad = \lim_n  p_{F}(u) \int_{X}G^*(x)
\biggl(\frac{\det\sigma_{F_{\rho}}(T_{\rho,\eta}x)}{\det\sigma
_F(x)}\biggr)^{1/2}
Y^{\eta}_{\rho}(x)[\det\sigma_F(x)]^{1/2} \nu_n^{F,\mbu}(dx)\\
&&\qquad = p_{F}(u) \int_{Z_{\mbu}}G^*(x)
\biggl(\frac{\det\sigma_{F_{\rho,\eta}}(T_{\rho,\eta}x)}{\det
\sigma_F(x)}\biggr)^{1/2}
Y^{\eta}_{\rho}(x) [\det\sigma_F(x)]^{1/2} \nu^{F,\mbu}(dx)\nonumber\\
&&\qquad = \int_{Z_{\mbu}}G^*(x)
\biggl(\frac{\det\sigma_{F_{\rho,\eta}}(T_{\rho,\eta}x)}{\det
\sigma_F(x)}\biggr)^{1/2}
Y^{\eta}_{\rho}(x) \,da^{Z_{\mbu}},\nonumber
\end{eqnarray}
which proves the following theorem.

%th4.6 #&#
\begin{theorem}
\label{theoremRNsurface}
Let $F \in D^{\infty-}_{1+\eps}$ satisfy the conditions of Theorem
\ref{theoremwatanabe} and
$\eta\in D^{\infty-}_{1+\eps}(X;H)$ be such that:
\begin{itemize}
\item$(I+\rho\eta)$ is one-to-one and onto when restricted to a
domain $B_{\eta}$ with complement
having $C^p_{\eps}$ capacity 0 for all $p$;
\item$(I_H+\rho\nabla\eta)$ is an invertible operator on $H$, when
restricted to $B_{\eta}$.
\end{itemize}
Then,
%
%e29 #&#
\begin{eqnarray}
\label{eqRNsurface}
\frac{da^{Z_{\mbu}}\circ T_{\rho,\eta}}{da^{Z_{\mbu}}}(x) & = &
\biggl(\frac{\det\sigma_{F_{\rho,\eta}}(T_{\rho,\eta}x)}{\det
\sigma_F(x)}\biggr)^{1/2} Y^{\eta}_{\rho}(x)
\nonumber
\\[-8pt]
\\[-8pt]
\nonumber
& \definedas& J^{Z_{\mbu}}_{\rho,\eta}.
\end{eqnarray}
\end{theorem}

%re4.7 #&#
\begin{remark}
\label{remarkconditions}
\textup{(1)} In order to better understand the above theorem, we shall now try
to simplify the expression
involved in \eqref{eqRNsurface}, for the simple case where $F=\delta
(h)$, for some $h\in H$,
and $\eta=\nabla F=h$, is a constant vector field. Clearly, $\nabla
\eta\equiv0$, implying
$DT_{\rho}=DT_{\rho}^{-1}=I_H$. Then the whole expression boils down to
\begin{eqnarray*}
J^{Z_{\mbu}}_{\rho,\eta} & = & Y^{\eta}_{\rho}(x)\\
& = & \exp\bigl(-\rho\delta(\eta)-\rho^2\|\eta\|^2_H/2\bigr).
\end{eqnarray*}

\begin{longlist}[(2)]
\item[(2)] The above expression in \eqref{eqRNsurface} can be rewritten as
\[
J^{Z_{\mbu}}_{\rho,\eta} = \biggl(
\frac{\det(\langle DT_{\rho,\eta}^{-1}(x)\nabla F_i(x),
DT_{\rho,\eta}^{-1}(x)\nabla F_j(x)\rangle_H)_{ij}}{\det(\langle
\nabla F_i(x),\nabla F_j(x)
\rangle_H)_{ij}}\biggr)^{1/2}Y^{\eta}_{\rho}(x),
\]
where $DT_{\rho}^{-1}$ is the operator given by $(I_H+\rho\nabla\eta)^{-1}$.

\item[(3)] From the definition of $Y^{\eta}_{\rho}(x)$, note that the
expression in
\eqref{eqRNsurface} is well defined as long as $\nabla\eta$ is a
Hilbert--Schmidt
class-valued operator acting on $H\times H$. Next,\vspace*{-1pt} in order to be able
to use the formula in
\eqref{eqRNsurface}, we need $J^{Z_{\mbu}}_{\rho,\eta}$ to be
integrable with
respect to the surface measure $da^{Z_{\mbu}}$.

\item[(4)] Note that the submanifold $Z_{\mbu}$ and the measure
$da^{Z_{\mbu}}$
are not dependent on~$F$. Therefore, to make the above
calculation simpler, we can choose an appropriate functional $F^{\prime}$,
such that
$\{F^{\prime}(x)=\bv\}=Z_{\mbu}$, for some $\bv\in\real^k$,
and that $\{\nabla F^{\prime}_i\}$ form an orthonormal
basis of the normal space of $Z_{\mbu}$.

\item[(5)] Finally, we note that $a^{Z_u} \circ T_{\rho,\eta}$ is defined
only up to capacity\vspace*{1pt} $C^{p}_{\eps}$ sets for any $p$. That is,
$C^p_{\eps}(A)=0$ implies
$a^{Z_u} \circ T_{\rho,\eta}(A)=0$. Hence, the image of the
discontinuities of $Z_u$ under
$T_{\rho,\eta}$ has $C^p_{\eps}$ capacity 0.
\end{longlist}
\end{remark}

%%-------------------------------------------------------------------------------------------------------
%%-------------------------------------------------------------------------------------------------------
%s4.3 #&#
\subsection{The set and its tube}
\label{subsecset}

Finally, we shall define the class of sets for which we shall prove a
tube formula, and, therefore,
define the GMFs. In our bid to keep the calculations much easier to
handle, we shall restrict our attention to
the unit codimensional case.

Continuing the way we have been defining subsets of the Wiener space
via Wiener functionals,
we shall start with a nondegenerate
Wiener functional $F\in D^{\infty-}_{2+\eps}(X;\real)$, such that
$F$ is an $H$-convex functional.
We shall write $A_{\mbu}= F^{-1}(-\infty,\mbu]$ for $\mbu\in\calO
$. This
is an $H$-convex set and its boundary $\partial A_{\mbu}$ is a smooth
unit codimensional submanifold of
the Wiener space.

For each $x \in A_{\mbu}$ define the support cone
\[
\support_x(A_{\mbu}) = \{h \in H\dvtx  \mbox{for any $\delta> 0,
 \ \exists  0 < \varepsilon< \delta$
such that $x+\varepsilon h \in A_{\mbu}$} \}
\]
and its dual, the (convex) normal cone
\[
N_x(A_{\mbu}) = \{h \in H\dvtx  \langle h, h' \rangle\leq0\
\forall  h' \in
\overline{\support_x(A_{\mbu})} \},
\]
where $\overline{\support_x(A_{\mbu})}$ is the closure of $\support
_x(A_{\mbu})$ in $H$.

The following lemmas describe some of the properties
of the tube around~$A_{\mbu}$. For clarity we state the results only
for the case
of unit codimension, though similar statements hold for
$k$-codimension.

Define the smooth points of $F$,
%
%e30 #&#
\begin{eqnarray}
\label{eqsmoothpoint}
\operatorname{Sm}(F) &\definedas& \Bigl\{x\dvtx  \grad F(x) \neq0\ \forall
h_1, h_2 \in H_{\infty}\
\exists  \varepsilon_n(h_1, h_2),
\downarrow0
\nonumber
\\[-8pt]
\\[-8pt]
\nonumber
&&\hspace*{16pt}\qquad \mbox{such that} \lim_{n \rightarrow\infty}
T_2(F,x,h_1,h_2,\varepsilon_{n}) = 0 \Bigr\},
\end{eqnarray}
where
\begin{eqnarray*}
T_2(F,x,h_1,h_2,\varepsilon) &=& F(x+h_1) + \varepsilon\langle\grad
F(x+h_1), h_2 \rangle_H
\\
&&{} + \varepsilon^2 \langle\grad^2 F(x+h_1), h_2 \otimes h_2
\rangle_{H \otimes H}
- F(x+h_1 + \varepsilon h_2)
\end{eqnarray*}
is the difference between the second-order Taylor expansion of
$F(x+h_1+\varepsilon h_2)$ evaluated
at $x+h_1$ and its true value.

%le4.8 #&#
\begin{lemma}
Suppose $F \in D^{\infty-}_{2+\eps}(X;\real)$.
Then, $C^p_{\eps}(\operatorname{Sm}(F)^c) = 0\ \forall p>1$.
At every $x \in\partial A_{\mbu} \cap\operatorname{Sm}(F)$,
%
%e31 #&#
\begin{equation}
\label{eqnormalspace}
N_x(A_{\mbu}) = \{ c \grad F(x)\dvtx  c \geq0 \}.
\end{equation}
\end{lemma}

\begin{pf} The first conclusion follows
essentially directly from Corol-\break lary~\ref{corsmoothness}.
Suppose now that $x \in\operatorname{Sm}(F)$. Then, for any $h \perp\grad
F(x), \|h\| \leq K$
we can find a sequence $H_{\infty} \ni h_n \overset{n\rightarrow
\infty}{\rightarrow} h$
satisfying $\langle h_n, \grad F(x) \rangle_H < -1/n$. Because the 2nd
order Taylor expansion
holds at $x$, we see that $h_n \in\support_x(A_{\mbu})$ for all~$n$.
Hence, $h \in\overline{\support_x(A_{\mbu})}$. This is enough to
conclude that any
$\eta\in N_x(A_{\mbu})$ is parallel to $\grad F(x)$. It is not hard
to see that it must
therefore be a positive multiple of $\grad F(x)$.
\end{pf}

%le4.9 #&#
\begin{lemma}
\label{lem121}
Suppose $F \in D^{\infty-}_{2+\eps}$ is $H$-convex. Then, for each $r
> 0$, the restriction of
\[
x \mapsto x + r \grad F(x) / \|\grad F(x)\| \definedas x + r \eta_x
\]
to $\operatorname{Sm}(F) \cap\partial A_{\mbu}$ is one-to-one in the sense that
for each $x \in\operatorname{Sm}(F) \cap\partial A_{\mbu}$
\[
\{y \in\operatorname{Sm}(F) \cap\partial A_{\mbu}\dvtx  \|y-(x + r \eta
_x)\|_{H} \leq r \} = \varnothing.
\]
\end{lemma}
\begin{pf} Given $x \in\operatorname{Sm}(F) \cap\partial A_{\mbu}$,
suppose such a $y$ exists with
$\|x+r \eta_x - y\| < r$. As $y$ is a smooth point of $F$, we can find
some $h \in H_{\infty}$
such that $\|x+r \eta_x - (y+h)\|_H < r$ and $F(y+h) < u$ with $F$
continuous at \mbox{$y+h$}.
Choose $\nu(x,y) \in H_{\infty}$ such that $x+\nu(x,y)$ is
arbitrarily close
to $y+h$. Then, by continuity of $F$ on $\operatorname{Sm}(F)$, $F(x+\nu(x,y))<u$
and $\|x+r \eta_x - (x+\nu(x,y))\|_H = \|r \eta_x - \nu(x,y)\|_H < r$.
Note that this implies $\langle\nu(x,y), \eta_x\rangle_H > 0$, or,
alternatively,
$\nu(x,y) \notin\support_x(A_{\mbu})$.\vadjust{\goodbreak}

Now consider the restriction of $F$ to the line segment joining
$[x,x+\nu(x,y)]$,
denoted by
\[
f(t) = F\bigl(x + t \bigl(x-\nu(x,y)\bigr)\bigr), \qquad 0 \leq t \leq1,
\]
which, by Remark~\ref{remlinearstructure}, is continuous,
twice-differentiable and convex on a
dense subset of $t \in[0,1],$ hence, we can find a continuous,
twice-differentiable convex function
$\tilde{f}$ on all of $[0,1]$ that agrees with $f$ on this dense subset.

There are two possibilities,
the first being that $\tilde{f}(t) \leq u$ for all $t \in[0,1]$. This
would imply
$\eta(x,y) \in\support_x(A_{\mbu})$, contradicting our previous observation.\vspace*{1pt}
The second alternative is that there exists $t$ such that $\tilde
{f}(t) > u$.
However, $\tilde{f}(0) = u, \tilde{f}(1) < u$ and this would violate
convexity.
By contradiction, there can be no such $y$.
This proves the assertion that there are no points $y$ of distance
strictly less than $r$ to $x+r \eta_x$. Now, suppose there exists a smooth
point $y \neq x$ of distance exactly $r$ from $x+r \eta_x$. Then, for
any $\delta> 0$
it is not hard to show that
\[
\bigl\|y - \bigl(x+(\delta+r)\eta_x\bigr)\bigr\|_H < \delta+ r,
\]
but we just proved that there can be no such $y$.
\end{pf}

We are now in a position to define the tube
\[
\Tube(A_{\mbu},\rho) = \{y \in X\dvtx  \exists x \in A_{\mbu}, \|
y-x\|_H \leq\rho\}
= \{y \in X\dvtx  d_H(y, A_{\mbu}) \leq\rho\},
\]
where the distance function is defined as
%
%e32 #&#
\begin{equation}
\label{eqdefdistance}
d_H(y, A_{\mbu}) = \inf_{h \in H_{\infty}\dvtx  y-h \in A_{\mbu}}\|h\|_H.
\end{equation}
The level sets of the distance function are hypersurfaces at distance $r$,
%
%e33 #&#
\begin{equation}
\label{eqlevelset}
\partial A_{\mbu}^r = \{y \in X\dvtx  d_H(y, A_{\mbu}) = r \}.
\end{equation}
Lemma~\ref{lem121} asserts that the restriction of $x \mapsto x + r
\eta_x$
to $A_{\mbu} \cap\operatorname{Sm}(F)$ is one-to-one. On the image of $\operatorname
{Sm}(F)$, its inverse is easily defined
as $x + r \eta_x \mapsto(x, \eta_x)$,
and, as noted in the remarks following Theorem~\ref{theoremRNsurface},
the image of $\partial A_{\mbu} \cap\operatorname{Sm}(F)$ has $C^p_{\eps
}$-capacity 0.
Hence, up to a set of $C^p_{\eps}$-capacity $0$, it is a bijection and
Theorem~\ref{theoremRNsurface} can be applied to study
the surface measure of~$\partial A_{\mbu}^r$.\looseness=-1

Moreover, the following theorem further corroborates the fact that the
change of measure formula established in
Theorem~\ref{theoremRNsurface} is the appropriate result to use in
order to obtain a tube formula, as will be
seen later.

%th4.10 #&#
\begin{theorem}
\label{theoremtubecapacity}
Let $C^{\infty-}_{\eps}(A)=0$, then under hypotheses \textup{(H2)} and \textup{(H3)} of~\cite{RenRockner05}, for $\eps_1<\eps$,
\[
C^{\infty-}_{\eps_1}\bigl(A\oplus B_H(0,r)\bigr)=0,
\]
where $B_H(0,r)$ is a ball in $H$ centered at $0$ with radius $r$.
\end{theorem}

Since capacities are continuous from below, it suffices to prove
that\break
$C^{\infty-}_{\eps_1}(A\oplus B_{E_n}(0,r))=0$,
for each $n$, whenever $C^{\infty-}_{\eps}(A)=0$, where
$B_{E_n}(0,r)$ is a ball of radius $r$, centered at $0$,\vadjust{\goodbreak}
in the vector space $E_n=\operatorname{span}(h_1,\ldots,h_n)$, where $\{h_i\}
_{i\ge1}$ is the orthonormal basis of $H$.
Also, note that the proof is given for an open subset $A$ of the Wiener
space $X$, but, using the arguments of
\cite{RenRockner05}, we can extend it to general subsets of the
Wiener space.

Before proving the above theorem, we shall, first, obtain some
estimates on functionals derived
from the Wiener functionals. Note that
\[
A\oplus B_{E_n}(0,r)=\bigl\{\bigl(A+\bigl\lla s,h^{(n)}\bigr\rra\bigr)\dvtx  s\in B_{\real
^n}(0,r)\bigr\},
\]
where $\lla s,h^{(n)}\rra=\sum_{i=1}^ns_ih_i$. Further, for the later
part, we shall denote
$I_n\subset B_{\real^n}(0,r)$ as the set of all rationals in the set
$B_{\real^n}(0,r)$.
The following result is, essentially, an extension of Theorem 2.1 of
\cite{RenRockner00}.

%th4.11 #&#
\begin{theorem}
\label{theoremHolder}
Let $f\in D^p_{\alpha}(X)$ for $\alpha\in(1/p,1)$, and $\real^n\ni
t\mapsto\xi(t,\cdot)=f(\cdot+\lla t,h^{(n)}\rra)$,
such that $|t|\le T$, for some fixed $T$, that is, $t$ belongs to some
large enough cube.
Then for all $p'\in(1/\alpha,p)$ there exists a $C=C(p,p',\alpha
,T)$, such that
\[
\|\xi(t)-\xi(s)\|_{p'}\le C \|f\|_{p,\alpha} |t-s|^{\alpha}.
\]
\end{theorem}

\begin{pf} Before we shall start proving the above result, we shall
recall that the estimates of Lemma 4.1 of
\cite{RenRockner00} remain unchanged in our setup. Now we need an
estimate analogous to the one obtained in Lemma
4.2 of~\cite{RenRockner00}, for which we recall the Ramer's change of
measure formula,
%
%e34 #&#
\begin{eqnarray}
\label{eqGlip}
&&\bigl\|G\bigl(\cdot+\bigl\lla t_2,h^{(n)}\bigr\rra\bigr)-G\bigl(\cdot+\bigl\lla
t_1,h^{(n)}\bigr\rra\bigr)\bigr\|_{p'}\nonumber\\[-2pt]
&&\qquad = \biggl\|G\biggl(\cdot+\frac{1}{2}\bigl\lla t_1+t_2,h^{(n)}\bigr\rra+\frac{1}{2}\bigl\lla
t_2-t_1,h^{(n)}\bigr\rra\biggr)\nonumber\\[-2pt]
&&\hspace*{3pt}\qquad\quad{}  -G\biggl(\cdot+\frac{1}{2}\bigl\lla t_1+t_2,h^{(n)}\bigr\rra-\frac
{1}{2}\bigl\lla t_2-t_1,h^{(n)}\bigl\rra\biggr)\biggr\|_{p'}\\[-2pt]
&&\qquad = \biggl(\int_X \biggl|G\biggl(x+\frac{1}{2}\bigl\lla t_2-t_1,h^{(n)}\bigr\rra\biggr)-G\biggl(x-\frac
{1}{2}\bigl\lla t_2-t_1,h^{(n)}\bigr\rra\biggr)\biggr|^{p'}\nonumber\\[-2pt]
& &\hspace*{50pt} {}\times \exp\biggl(-\frac{1}{2}\biggl\|\frac{1}{2}\bigl\lla t_1+t_2,h^{(n)}\bigr\rra\biggr\|
_H^2-\delta\biggl(\frac{1}{2}\bigl\lla t_1+t_2,h^{(n)}\bigr\rra\biggr)\biggr)
 \mu(dx)\biggr)^{1/p'}.\nonumber
\end{eqnarray}
Now writing $h_1=|t_2-t_1|^{-1}\lla t_2-t_1,h^{(n)}\rra$,
$h_2=|t_1+t_2|^{-1}\lla t_1+t_2,h^{(n)}\rra$, and
$Y^{h_2}_{|t_1+t_2|/2}=\exp[-\||t_1+t_2|h_2/8\|_H^2-\delta
(|t_1+t_2|h_2/2)]$, we
can rewrite the above as
%
%e35 #&#
\begin{eqnarray}
&&\bigl\|G\bigl(\cdot+\bigl\lla t_2,h^{(n)}\bigr\rra\bigr)-G\bigl(\cdot+\bigl\lla
t_1,h^{(n)}\bigl\rra\bigr)\bigr\|_{p'}\nonumber\\[-2pt]
&&\qquad = \biggl(\int_X \biggl|G\biggl(x+\frac{1}{2}|t_2-t_1|h_1\biggr)-G\biggl(x-\frac
{1}{2}|t_2-t_1|h_1\biggr)\biggr|^{p'}\\[-2pt]
&&\hspace*{151pt}\qquad{}\times Y^{h_2}_{|t_1+t_2|/2}(x)
\mu(dx)\biggr)^{1/p'}.\nonumber
\end{eqnarray}
This reduces the above expression to the case dealt in \cite
{RenRockner00}. Therefore,
using the rest of the calculations of Lemma 4.2 of \cite
{RenRockner00}, and writing $G=T_af$,
where $\{T_a\}_{a\ge0}$ is the semigroup associated with the
Ornstein--Uhlenbeck operator~$L$,
we get the desired estimate expressed as
\[
\bigl\|T_af\bigl(\cdot+\bigl\lla t_2,h^{(n)}\bigr\rra\bigr)-T_af\bigl(\cdot+\bigl\lla t_1,h^{(n)}\bigr\rra
\bigr)\bigr\|_{p'}
\le C(p,p',\alpha)\|T_af\|_{p,1} |t_2-t_1|.
\]
Thereafter, we can mimic the proof of Theorem 2.1 of \cite
{RenRockner00} and get the desired estimate.
\end{pf}

Now coming back to our case, let $e_A$ be the potential equilibrium of
$A$ (cf.~\cite{Shigekawa94}) and
$e_A\in D^{\infty-}_{\eps}$, therefore, there exists a $v_A\in
L^{\infty-}$, such that
\[
e_A=(I-L)^{-\eps/2}v_A\definedas(I-L)^{-\eps_1/2}v_{(\eps-\eps_1),A},
\]
where $\eps_1$ is some number strictly smaller than $\eps$ and $L$ is
the Ornstein--Uhlenbeck operator.
Then, clearly,
\[
e_A\bigl(\cdot+\bigl\lla t,h^{(n)}\bigr\rra\bigr)=(I-L)^{-\eps/2}v_A\bigl(\cdot+\bigl\lla
t,h^{(n)}\bigr\rra\bigr)
\definedas(I-L)^{-\eps_1/2}v_{(\eps-\eps_1),A}\bigl(\cdot+\bigl\lla
t,h^{(n)}\bigr\rra\bigr).
\]
Now, writing $\xi(t)\definedas e_A(\cdot+\lla t,h^{(n)}\rra)$ and
$\xi_{(\eps-\eps_1)}(t)\definedas v_{(\eps-\eps_1),A}(\cdot+\lla
t,h^{(n)}\rra)$,
and also, in the process, choosing the appropriate quasi-continuous
redefinitions of the processes
$\xi$ and $\xi_{(\eps-\eps_1)}$, and choosing a large $p'$
(conditions on $p'$ will appear later),
such that by Kree--Meyer inequalities, we have
%
%e36 #&#
\begin{equation}
\label{eqxibound1}
\|\xi(t)-\xi(s)\|_{p',\eps_1} \le C\bigl\|\xi_{(\eps-\eps_1)}(t)-\xi
_{(\eps-\eps_1)}(s)\bigr\|_{p'}.
\end{equation}
Now using the above theorem with $f$ replaced by $v_{(\eps-\eps
_1),A}$, we get
%
%e37 #&#
\begin{eqnarray}
\label{eqxibound2}
&&\bigl\|\xi_{(\eps-\eps_1)}(t)-\xi_{(\eps-\eps_1)}(s)\bigr\|_{p'}\nonumber\\
 &&\qquad =
\bigl\|v_{(\eps-\eps_1),A}\bigl(\cdot+\bigl\lla t,h^{(n)}\bigr\rra\bigr)-v_{(\eps-\eps
_1),A}\bigl(\cdot+\bigl\lla t,h^{(n)}\bigr\rra\bigr)\bigr\|_{p'}
\nonumber
\\[-8pt]
\\[-8pt]
\nonumber
&&\qquad \le C\bigl \|v_{(\eps-\eps_1),A}\bigr\|_{p,(\eps-\eps_1)} |t-s|^{(\eps
-\eps_1)}\nonumber\\
&&\qquad =  C \|e_{A}\|_{p,\eps}  |t-s|^{(\eps-\eps_1)},\nonumber
\end{eqnarray}
where $p'\in(2/\eps,p)$.
Combining \eqref{eqxibound1} and \eqref{eqxibound2}, we get
%
%e38 #&#
\begin{equation}
\label{eqxibound3}
\|\xi(t)-\xi(s)\|_{p',\eps_1}\le C\|e_{A}\|_{p,\eps} |t-s|^{(\eps
-\eps_1)},
\end{equation}
which can be rewritten as
%
%e39 #&#
\begin{equation}
\label{eqShigekawacondition}
\sup_{s\neq t}\frac{\|\xi(t)-\xi(s)\|_{p',\eps
_1}^{p'}}{|t-s|^{p'(\eps-\eps_1)}} \le C\|e_A\|_{p,\eps}^{p'}.
\end{equation}
Now we can list the assumptions on the various indices as follows: we
start with any fixed $\eps_1<\eps$,
then choose a large enough $p$ such that $(\eps-\eps_1)\in(1/p,1)$,
and then we choose $p'$ such that
$p'\in(1/(\eps-\eps_1),p)$ and \mbox{$p'(\eps-\eps_1)>n$}. This can be
achieved by choosing $p$ and $p'$ of the
order of~$n$, in particular, choosing $p=an/(\eps-\eps_1)$ and
$p'=bn/(\eps-\eps_1)$, for
$a>b>0$ will do.
Then, using Theorem~3.4 of~\cite{Shigekawa94}, we get
%
%e40 #&#
\begin{equation}
\label{eqcapbound1}
C^p_{\eps_1}\Bigl(\sup_{s\neq t}|\xi(t)-\xi(s)|\Bigr) \le C \|
e_{A}\|_{p,\eps}^{p'}.
\end{equation}
\begin{pf*}{Proof of Theorem \protect\ref{theoremtubecapacity}} Now let us consider
%
%e41 #&#
\begin{eqnarray}
\label{eqcapbound2}
C^p_{\eps_1}\bigl(A\oplus B_{E_n}(0,r)\bigr)
& = & C^p_{\eps_1}\Bigl(\sup_{s\in B_{\real^n}(0,r)}1_A\bigl(\cdot+\bigl\lla
s,h^{(n)}\bigr\rra\bigr)\Bigr)\nonumber\\
& \le& C^p_{\eps_1}\Bigl(\sup_{s\in I_n}e_A\bigl(\cdot+\bigl\lla
s,h^{(n)}\bigr\rra\bigr)\Bigr)\qquad \mbox{as $e_A\ge1_A$}
\nonumber
\\[-8pt]
\\[-8pt]
\nonumber
& \le& C^p_{\eps_1}\Bigl(\sup_{s\in B_{\real^n}(0,r)}\xi(s)
\Bigr)\\
& \le& C^p_{\eps_1}\Bigl(\sup_{s\in B_{\real^n}(0,r)}|\xi(s)-\xi
(0)|+|\xi(0)|\Bigr).\nonumber
\end{eqnarray}
Now using \eqref{eqcapbound1}, we shall get
\[
C^p_{\eps_1}\bigl(A\oplus B_{E_n}(0,r)\bigr) \le(C+1) \|e_A\|_{p,\eps}^{p'}
= (C+1) (C^p_{\eps}(A))^{{p'}/{p}},
\]
which proves that $C^{\infty-}_{\eps_1}(A\oplus B_{E_n}(0,r))=0$, for
all $p>n(\eps-\eps_1)^{-1}$.
Now by the definition of the capacities and the hierarchy of the
Sobolev spaces, we shall have
$C^{\infty-}_{\eps_1}(A\oplus B_{E_n}(0,r))=0$, thereby proving the result.
\end{pf*}

Using the definition of the smooth points $\operatorname{Sm}(F)$ and $\Tube
(A_{\mbu},\rho)$, we can conclude
that
%
%e42 #&#
\begin{eqnarray}
\Tube(A_{\mbu},\rho) &=& \bigl[\bigl(A_{\mbu}\cap\operatorname{Sm}(F)\bigr)\oplus
B_H(0,\rho)\bigr]
\nonumber
\\[-8pt]
\\[-8pt]
\nonumber
&&{}\cup\bigl[\bigl(A_{\mbu}\cap\{\operatorname{Sm}(F)\}^c\bigr)\oplus B_H(0,\rho)\bigr].
\end{eqnarray}
Using the above calculations, we have
%
%e43 #&#
\begin{equation}
\label{eqtubeequivalence}
\mu(\Tube(A_{\mbu},\rho))=\mu\bigl(\bigl(A_{\mbu}\cap\operatorname
{Sm}(F)\bigr)\oplus B_H(0,\rho)\bigr),
\end{equation}
since $C^p_{\eps_1}((A_{\mbu}\cap\{\operatorname{Sm}(F)\}^c)\oplus
B_H(0,\rho))=0$, implying that
the $\mu$-measure of the set is zero. Therefore, it is enough, for the
tube formula, to consider the
set $((A_{\mbu}\cap\operatorname{Sm}(F))\oplus B_H(0,\rho))$, on which the
transformation $x\mapsto x+\eta_x$
is well defined up to $C^p_{\eps}$-zero sets, and, hence, we can use
the change of measure formula for the surface
areas given in Theorem~\ref{theoremRNsurface}.

%%------------------------------------------------------------------------------------------------------------------------------------------------
%%------------------------------------------------------------------------------------------------------------------------------------------------

%s5 #&#
\section{A Wiener tube formula}
\label{sectubeformula}

After setting up the basics, definitions and the conditions, concerning
a tube formula
in the Wiener space, we shall finally prove one of the main results of
this paper, which can be stated
in the form of the following theorem.\vadjust{\goodbreak}

%th5.1 #&#
\begin{theorem}
\label{theoremmain}
Let $F\in D^{\infty-}_{2+\delta}(X;\real)$ be an $H$-convex Wiener
functional such that it
satisfies all the regularity conditions of Theorem \ref
{theoremexistence}, and
$A_{\mbu}=F^{-1}(-\infty,\mbu]$, then
\[
\mu(\Tube(A_{\mbu},\rho))=\minkmu_0(A_{\mbu})+\sum_{j=1}^{\infty
}\frac{\rho^j}{j!}\minkmu_j(A_{\mbu}),
\]
where $\minkmu_j(A_{\mbu})$ are the infinite dimensional versions of
Gaussian Minkowski
functionals and, as usual, $\minkmu_0(A_{\mbu})=\mu(A_{\mbu})$.
\end{theorem}

\begin{pf} Let us start with recalling the definition of the {\it
outward} pointing normal space $N_x(A_{\mbu})$ from
\eqref{eqnormalspace} and
writing $N(A_{\mbu})=\bigcup_{x\in A_{\mbu}}N_x(A_{\mbu})$. Then, let
us define a distance function,
$d_{A_{\mbu}}\dvtx \Tube(A_{\mbu},\rho)\to\real$, such that
for $x\in\Tube(A_{\mbu},\rho)$ writing the ``residual'' as
\[
\hat{r}_x=
\mathop{\argmin}_{r\in\real; \eta\in N(A_{\mbu})}
d(x-r\eta,A_{\mbu}),
\]
the distance function $d_{A_{\mbu}}$
is given by $d_{A_{\mbu}}(x)= \|\hat{r}_x\|.$
Clearly, from the above definition, $d_{A_{\mbu}}^{-1}(0)=A_{\mbu}$.
Also, we can further express
$\Tube(A_{\mbu},\rho)$ as the disjoint union of $A_{\mbu}$ and
$\Tube^{+}(\partial A_{\mbu},\rho)$,
where $\Tube^{+}(\partial A_{\mbu},\rho)=\Tube(\partial A_{\mbu
},\rho)\cap A_{\mbu}^c$. Thus,
%
%e44 #&#
\begin{equation}
\label{eqtubedisjointsum}
\mu(\Tube(A_{\mbu},\rho)) = \mu(A_{\mbu}) + \mu(\Tube
^{+}(\partial A_{\mbu},\rho)).
\end{equation}
Now using the Wiener space version of Federer's co-area formula as it
appears in~\cite{AiraultMalliavin}, we shall obtain
%
%e45 #&#
\begin{equation}
\label{eqvolumetubelayers}
\mu(\Tube^{+}(\partial A_{\mbu},\rho))
= \int_0^{\rho}\int_{d_{A_{\mbu}}^{-1}(r)}
(\sigma_{d_{A_{\mbu}}}(x))^{-1} \,da^{\partial^{+} A_{\mbu}^r} \,dr,
\end{equation}
where\vspace*{1pt} $\partial^{+} A_{\mbu}^r=d_{A_{\mbu}}^{-1}(r)\cap A_{\mbu}^c$
are the level sets of the distance
function $d_{A_{\mbu}}$ in the outward direction.
Now note that $\nabla d_{A_{\mbu}}=\eta$, hence,
$\sigma_{d_{A_{\mbu}}}(x)=1$. Then let us define the
transformation $T_{r,\eta}\dvtx X\to X$, such that its restriction to
$A_{\mbu}$ is given by
$T_{r,\eta}(x)=x+r\eta.$
Clearly, $T_{r,\eta}(\partial A_{\mbu})=\partial A_{\mbu}^r$.
Then, we shall use our change of measure formula for surfaces on $\int
_{d_{A_{\mbu}}^{-1}(r)}$
to further simplify the expression in \eqref{eqvolumetubelayers} to obtain
\[
\mu(\Tube^{+}(\partial A_{\mbu},\rho)) = \int_0^{\rho}\int
_{A_{\mbu}} J^{\partial A_{\mbu}}_{r,\eta}
 \,da^{\partial A_{\mbu}} \,dr = \int_0^{\rho}\int_{A_{\mbu}}
Y^{\eta}_r  \,da^{\partial A_{\mbu}} \,dr,
\]
where terms $J^{\partial A_{\mbu}}_{r,\eta}$ and $Y^{\eta}_r$ are as
they appear in Theorem~\ref{theoremRNsurface}.

Now using a Taylor series expansion for $Y^{\eta}_r$ with respect to
$r$, we can rewrite
the above expression as
\begin{eqnarray*}
&&
\mu(\Tube^{+}(\partial A_{\mbu},\rho))\\
&&\qquad = \sum_{j=0}^{\infty}\frac{\rho^{j+1}}{(j+1)!}\int_{A_{\mbu}}
\frac{d^j}{dr^j}\bigl(\det_2(I_H+r\nabla\eta)
\exp\bigl(-r\delta(\eta)-r^2/2\bigr)\bigr)\bigg|_{r=0} \,da^{\partial A_{\mbu}}.
\end{eqnarray*}
We note here that $\rho$ must be within the radius of convergence of
the Taylor series of
$Y^{\eta}_r$,
which in turn will ensure the convergence of the above series.

Finally, plugging the above expression in \eqref
{eqtubedisjointsum}, we get
%
%e46 #&#
\begin{eqnarray}
&&\mu(\Tube(A_{\mbu},\rho))\nonumber\\
&&\qquad = \mu(A_{\mbu})
\nonumber
\\[-8pt]
\\[-8pt]
\nonumber
&&\qquad\quad{}+ \sum_{j=1}^{\infty}\frac{\rho^{j}}{j!}\int
_{A_{\mbu}}
\frac{d^j}{dr^j}\bigl(\det_2(I_H+r\nabla\eta)
\exp\bigl(-r\delta(\eta)-r^2/2\bigr)\bigr)\bigg|_{r=0} \,da^{\partial A_{\mbu
}}\\
&&\qquad = \mu(A_{\mbu})+ \sum_{n=1}^{\infty}\frac{\rho^n}{n!}\minkmu
_n(A_{\mbu}),\nonumber
\end{eqnarray}
where $\minkmu_n(A_{\mbu})$ are Gaussian Minkowski functionals of the
infinite dimensional
set $A_{\mbu}$, given by
%
%e47 #&#
\begin{equation}
\label{eqdefGMF}
\qquad\minkmu_n(A_{\mbu}) = \int_{A_{\mbu}}
\frac{d^n}{dr^n}\bigl(\det_2(I_H+r\nabla\eta)
\exp\bigl(-r\delta(\eta)-r^2/2\bigr)\bigr)\bigg|_{r=0} \,da^{\partial A_{\mbu}},
\end{equation}
which proves the theorem.
\end{pf}

%%---------------------------------------------------------------------------------------------------------------------
%%---------------------------------------------------------------------------------------------------------------------

%s6 #&#
\section{Applications}
\label{secapp}

In this section we shall invoke the existential results from the
previous section to
obtain a kinematic fundamental formula akin to the one obtained in
Theorem $15.9.5$ of
\cite{RFG}, though, for a larger class of random fields.

Let us consider a real-valued random field $f$ defined on a compact
Riemannian manifold $M$
equipped with a metric $\tau$. Then the modulus of continuity $\Xi$
of a function $F\dvtx M\to\real$
is defined as
\[
\Xi_F(\eta)\definedas\sup_{\tau(x,y)\le\eta}|F(x)-F(y)|
\]
for all $\eta>0$.
Continuing the setup introduced in the example stated in Section \ref
{secintro}, we shall consider
a specific class of random fields $f$ which can be represented as
%
%e48 #&#
\begin{equation}
\label{eqrandomfield}
f(x)=\sum_{i=1}^N\int_0^1 V_i(B_i^x(s))  \,dB_i^x(s),
\end{equation}
where the integral\vspace*{1pt} is to be interpreted in the It\^o sense, each
$V_i\dvtx \real\to\real$
is a smooth function, and $B^x(\cdot)=(B^x_i(\cdot))_{i=1}^N$
is a $\real^N$-valued, zero-mean Gaussian process defined on $M\times
[0,1]$, whose covariance is given by
%
%e49 #&#
\begin{equation}
\label{eqfcov}
E(B_i^x(s)B_j^y(t))=(s\wedge t) C_{ij}(x,y) = (s\wedge t) C(x,y),
\end{equation}
where $C\dvtx M\times M\to\real$ is a smooth function, such that for each
fixed $t\in[0,1]$, the field
$B^{\cdot}(t)$ is an isotropic
Gaussian field over $M$ (see Sections 5.7 and~5.8 of~\cite{RFG}).\vadjust{\goodbreak}

The following are the basic assumptions on the functions $V_i$ and the
Gaussian process $B^x(t)$:
\begin{longlist}[(A1)]
\item[(A1)] $V_i\in C^4, 1\le i\le N$, the class of all
$4$-continuously differentiable functions.
\item[(A2)] Writing $V^{(k)}$ as the $k$th derivative of $V$, let us define
\[
C_i(k,s,x,y)\definedas\sup_{0\le\a\le1} V^{(k)}_i\bigl(\a
B^x(s)+(1-\a)B^y(s)\bigr)
\]
for any $x,y\in M$
and $k=0,1,2,3,4$. Then, for some $p\gg\dim(M)$ and for all
$i=1,\ldots,N$,
\[
\sup_{0\le s\le1}\|C_i(k,s,x,y)\|_p^p = c_{i,k}(x,y,p)<\infty.
\]
Also, $\sup_{x\neq y}c_{i,k}(x,y,p)< \infty,$ for all $k=0,1,2,3,4$.
Note that
this is satisfied whenever the $V_i$'s are $C^4$ with polynomial growth.
\item[(A3)] For each $r\ge1$, there exists a constant $m_r$, such that
\[
E|B^x(s)-B^y(s)|^r \le m_r|x-y|^r \qquad \forall0\le s\le1,
\]
where $m_r$ depends solely on $r$.
\item[(A4)] All the above assumptions also hold true with $B^x(s)$
replaced by $\nabla B^x(s)$ and
$\nabla^2 B^x(s),$ respectively.
\end{longlist}

Let us define the \textit{excursion set} $A_u$ corresponding to $f\dvtx M\to
\real$ as
\[
A_u(f;M)\definedas\{x\in M\dvtx  f(x)\ge u\}.
\]
Also, note that writing $F(\omega) = \int_0^1 V(\omega_s) \,d\omega
_s,$ one can consider the random field
defined above as $f(x) = F(B^x)$.

%th6.1 #&#
\begin{theorem}
\label{theoremGKF}
Let $M$ be a $m$-dimensional manifold and $f$ be a random field defined
on $M$, represented as in
\eqref{eqrandomfield}, and satisfying the conditions \textup{(A1)--(A4)}.
Also, let $f(x)$ and $\nabla f(x)$
be nondegenerate in the sense of Malliavin, for some $x\in M$,
and that the corresponding Wiener functional $F$ satisfies the
exponential moment condition specified in Theorem
\ref{theoremexistence}. Then writing $A_u(f;M)$ as the excursion set
for the random field $f$, and
$\lips_i(\cdot)$ as the $i$th Lipschitz--Killing curvature under the
Gaussian induced metric, we have
%
%e50 #&#
\begin{equation}
\label{eqGKF}
E(\lips_0(A_u(f;M)))=\sum_{j=0}^{m}(2\pi)^{-j/2} \lips_{j}(M)
\minkmu_j(F^{-1}[u,\infty)),
\end{equation}
where $F^{-1}([u,\infty))$ is a subset of the Wiener space $X$ with
$\minkmu_j$ as its
GMF, as defined in the previous section.
\end{theorem}

We shall approximate the LHS of \eqref{eqGKF} such that it has the
form of the RHS and that this approximation
of the RHS indeed converges to the RHS of \eqref{eqGKF}.\vadjust{\goodbreak}

At this point, we note that all the results obtained below are for the
case of $N=1$, whereas using
similar methods, the same results are true for general $N$.
We shall start proving Theorem~\ref{theoremGKF} by first listing some
regularity properties of
field~$f$ defined in \eqref{eqrandomfield}, in the form of the
following theorem.

%th6.2 #&#
\begin{theorem}
\label{theoremfieldproperties}
Let the random field $f$ be as defined in \eqref{eqrandomfield},
such that it also satisfies \textup{(A1)--(A4)}, then:
\begin{longlist}[(a)]
\item[(a)] $F\in D^{\infty-}_{3}(X;\real)$, and under the assumption
of nondegeneracy of $F$, the density $p_F$ of
$F$ is bounded,
\item[(b)] $f$ is continuous, and that for any $\eps>0$
\[
P\bigl(\Xi_f(\eta)>\eps\bigr) = o\bigl(\eta^{\dim(M)}\bigr) \qquad   \mbox{as }  \eta
\downarrow0.
\]
Also, the same is true for $\nabla f$ and $\nabla^2 f$.
\end{longlist}
\end{theorem}

\begin{pf} Clearly, for the Wiener functional $F=\int_0^1V(B(t))
\,dB(t)$, we have $DF\in H$.
Therefore, by definition,\vspace*{-2pt} there exists a unique $(\widehat{DF})\in
L^2([0,1])$ such that
$D_rF \definedas(DF)(r)=\int_0^r (\widehat{D_sF}) \,ds$.
Using this notation, we shall have
\begin{eqnarray*}
(\widehat{D_rF})& =& V(B(r)) + \int_0^1 V^{(1)}(B(t))1_{[0,t]}(r) \,dB(t),
\\
(\widehat{D_s\widehat{(D_rF)}}) & = &
V^{(1)}(B(r))1_{[0,r]}(s) + V^{(1)}(B(s))1_{[0,s]}(r)\\
& &{} + \int_0^1 V^{(2)}(B(t))1_{[0,t]}(r)1_{[0,t]}(s) \,dB(t).
\end{eqnarray*}
Clearly, due to the moment conditions imposed on $V$ and its
derivatives, we can conclude that
$F\in D^{\infty-}_3(X;\real)$, and the boundedness of the density
$p_F$ follows using Proposition $2.1.1$
of~\cite{Nualart-book}.

Now to prove continuity of $f$ and its derivatives, we shall use
Kolmogorov's continuity criterion for processes
defined on smooth Riemannian manifolds.
Note that Kolmogorov's continuity criterion is usually stated for
processes with Euclidean
parameter space, but
since continuity is a local phenomena, thus, it can easily be extended
to processes defined
on smooth locally Euclidean spaces.
Therefore, we present the proof of continuity related results for the
field $f\circ\phi^{-1}$ where~$\phi$ is the local chart,
but we shall suppress the chart map, and will write $f$ for both the
field $f$
and its counterpart $f\circ\phi^{-1}$.

In order to use Kolmogorov's continuity theorem, we must obtain $L^p$
estimates for $(f(x)-f(y))$.
Writing $V^*$ as any antiderivative of $V$, we have
\[
V^*(B^x(1))=V^*(B^x(0))+\int_0^1 V(B^x(s)) \,dB^x(s) + \int_0^1
V^{(1)}(B^x(s)) \,ds.
\]
Thus, for $p\ge1$, there exist $m_{1,p}$ and $m_{2,p}$ such that
\begin{eqnarray*}
\|f(x)-f(y)\|_p^p &= &E|f(x)-f(y)|^p\\
& \le& m_{1,p} E|V^*(B^x(1))-V^*(B^y(1))|^p
\\
&&{}+ m_{2,p} E\biggl|\int_0^1 \bigl(V^{(1)}(B^x(s))- V^{(1)}(B^y(s))
\bigr) \,ds\biggr|^p\\
& \le& m_{1,p} E\Bigl|\sup_{\a}V[\a B^x(1)+(1-\a)B^x(1)
]\times\bigl(B^x(1)-B^y(1)\bigr)\Bigr|^p\\
&&  {} + m_{2,p} \int_0^1 E\bigl|V^{(1)}(B^x(s))- V^{(1)}(B^y(s))\bigr|^p
\,ds\\
& \le &m_{1,p} E\bigl|C(0,1,x,y)\bigl(B^x(1)-B^y(1)\bigr)\bigr|^p
\\
&&{}+ m_{2,p} E\bigl\|C(2,\cdot,x,y)\bigl(B^x(\cdot)- B^y(\cdot)\bigr)\bigr\|^p_{L^p[0,1]}.
\end{eqnarray*}
Now choosing $p_1,p_2 >0$ such that $p_1^{-1}+p_2^{-1}=p^{-1}$, we get
\begin{eqnarray*}
\|f(x)-f(y)\|_p^p & \le& m_{1,p}\bigl (\|C(0,1,x,y)\|_{p_1}\|
B^x(1)-B^y(1)\|_{p_2}\bigr)^p\\
& &  {} + m_{2,p} \int_0^1 \bigl(\|C(2,s,x,y)\|_{p_1}\|
B^x(s)-B^y(s)\|_{p_2}\bigr)^p \,ds\\
& \le& m_{1,p} c_0(x,y,p_1)^{p/p_1} m_{p_2}^{p/p_2}|x-y|^p\\
& &  {} + m_{2,p} c_2(x,y,p_1)^{p/p_1} m_{p_2}^{p/p_2}|x-y|^p.
\end{eqnarray*}
Next, fixing
\[
M(p,p_1,p_2) = \sup_{x\neq y}\bigl(m_{1,p} c_0(x,y,p_1)^{p/p_1} m_{p_2}^{p/p_2}
+ m_{2,p} c_2(x,y,p_1)^{p/p_1} m_{p_2}^{p/p_2}\bigr),
\]
we have
%
%e51 #&#
\begin{equation}
\|f(x)-f(y)\|\le M(p,p_1,p_2)|x-y|^p.
\end{equation}
For large enough $p$ we can use Theorem $1.4.1$ in~\cite{Kun97} to
deduce that there exists $\tilde{f}$,
which is the continuous modification of $f$. Abusing the notation, we
shall write $f$ for $\tilde{f}$.
Also, using the same result, we can infer that for any $\eps>0$, the
modulus of continuity of $f$ satisfies
\[
P\bigl(\Xi_f(\eta)>\eps\bigr) = o\bigl(\eta^{\dim(M)}\bigr) \qquad   \mbox{as }  \eta
\downarrow0.
\]

Note that we needed supremum of $c_2(x,y)$ to be bounded to prove the
continuity of $f$ and to control
its modulus of continuity. We can further conclude that the conditions
stated in \textup{(A1)--(A4)} suffice
to obtain similar results for the modulus of continuity of $\nabla f$
and $\nabla^2f$.
\end{pf}

Recall from~\cite{RFG} that $\lips_0$, also known as the
Euler--Poincar\'e characteristic, of the excursion set
$A_u(f;M)$ can be expressed as
\begin{eqnarray*}
\lips_0(A_u(f;M)) & = & \sum_{k=0}^m(-1)^k\#\{x\in M\dvtx f(x)\ge u,
\nabla f(x)=0, \operatorname{index}(\nabla^2 f)=k\}\\
& = & \sum_{k=0}^m (-1)^k \mu_k.
\end{eqnarray*}

Now using the \textit{expectation metatheorem} (Theorem $11.2.1$ of \cite
{RFG}), and replacing~$G$ and $H$ by $\nabla f$ and $(\nabla^2f,f)$,
respectively, and $B$ by $D_k\times[u,\infty)$,
where~$D_k$ is the space of $m\times m$ matrices with index $k$, we can
obtain a formula for the expected
value of $\mu_k$ as defined above.
However, in order to use this result for our purpose, we must also
check the conditions involving
conditional densities, for which we refer to Theorem $4.1$ of \cite
{NualartZakai}.
Thus, using these results, we can write
%
%e52 #&#
\begin{eqnarray}
\label{eqexpectationmeta}
&&
E(\lips_0(A_u(f;M)))
\nonumber
\\[-8pt]
\\[-8pt]
\nonumber
& &\qquad= \int_M E\bigl(\det(-\nabla^2 f(x)) 1_{[u,\infty)}(f(x))|\nabla
f(x)=0\bigr)
 p_{\nabla f(x)}(0) \,dx.
\end{eqnarray}

Next, in order to construct an approximating sequence to the LHS of
\eqref{eqGKF} and appeal to the
results in~\cite{RFG}, we shall use a cylindrical approximation of
$f(x)$. Let $\{(i/n,(i+1)/n]\}_{i=0}^{n-1}$
be a partition of $(0,1]$, then define
\[
f_n(x) = \sum_{i=0}^{n-1} V\bigl(B^x(i/n)\bigr)\bigl(B^x\bigl((i+1)/n\bigr)-B^x(i/n)\bigr).
\]
Standard results from stochastic analysis ensure the convergence of
$f_n(x)$ to $f(x)$.
Moreover, note that $(B^x((i+1)/n)-B^x(i/n))_{i=0}^{n-1}$ forms an i.i.d.
$0\le i\le(n-1)$. Therefore, we can write
\[
f_n(x)=F_n\bigl(y^{(n)}_1(x),\ldots,y^{(n)}_n(x)\bigr),
\]
where $y^{(n)}_{i+1}(x)$ are i.i.d. with the same distribution as
$\sqrt{n}(B^x((i+1)/n)-B^x(i/n))$
and $F_n$ is the appropriately defined real-valued function.
Under the conditions imposed on $f$ for the expectation metatheorem to
be true,
$f_n$ also becomes a valid candidate to apply the metatheorem, thereby
giving us
%
%e53 #&#
\begin{eqnarray}
\label{eqexpectationmetan}
\qquad &&E(\lips_0(A_u(f_n;M)))
\nonumber
\\[-8pt]
\\[-8pt]
\nonumber
&&\qquad =\sum_{k=0}^N\int_M E\bigl(\det(-\nabla^2 f_n(x)) 1_{[u,\infty
)}(f(x))|\nabla f_n(x)=0\bigr)
 p_{\nabla f_n(x)}(0) \,dx.
\end{eqnarray}
Using Theorem 15.9.5 of~\cite{RFG} for the random field $f_n$, we
shall have
%
%e54 #&#
\begin{equation}
\label{eqGKFn}
E(\lips_0(A_u(f_n;M)))=\sum_{j=0}^{m}(2\pi)^{-j/2} \lips_{j}(M)
\minkmu_j(F_n^{-1}[u,\infty)),
\end{equation}
where $F_n^{-1}[u,\infty)$ is a subset of $\real^n$.

%th6.3 #&#
\begin{theorem}
\label{theoremGMFconvergence}
Let $\{G_n\}_{n\ge1}$ be a sequence of real-valued Wiener functionals,
such that $G_n$ belongs to
the $n$th Wiener chaos,
and $G_n\to G$ in $D^{\infty-}_{3}$, for some $G\in D^{\infty-}_{3}$.
Also, let that each $G_n$ and $G$ satisfy all the assumptions of
Theorem~\ref{theoremexistence}, then
$\minkgk_j(G_n^{-1}[u,\infty))\to\minkmu_j(G^{-1}[u,\infty))$, as
$n\to\infty$.
\end{theorem}

\begin{pf} Using the definition of GMFs in Theorem \ref
{theoremexistence} and the convergence of the densities
$p_{G_n}$ to $p_G$, it suffices to prove that
\begin{eqnarray*}
&&E^{G_n=u}\bigl([\det(\sigma_{G_n})]^{1/2}\det_2(I_H+rD\eta
_n)\exp\bigl(-r\delta(\eta_n)-\tfrac{1}{2}r^2\bigr)\bigr)\\
&&\qquad \to E^{G=u}\bigl([\det(\sigma_G)]^{1/2}\det_2(I_H+rD\eta)\exp
\bigl(-r\delta(\eta)-\tfrac{1}{2}r^2\bigr)\bigr),
\end{eqnarray*}
where $\eta=DG/\|DG\|_H$ and $\eta_n=DG_n/\|DG_n\|$.
Writing
\[
A_n = \bigl([\det(\sigma_{G_n})]^{1/2}\det_2(I_H+rD\eta_n)\exp
\bigl(-r\delta(\eta_n)-\tfrac{1}{2}r^2\bigr)\bigr),
\]
and similarly defining $A$, we get
%
%e55 #&#
\begin{eqnarray}
\label{eqGMFcgs}
&&|E^{G_n=u}A_n - E^{G=u}A|\nonumber\\
&&\qquad = |E(A_n  \delta_u\circ G_n)-E(A  \delta
_u\circ G)|
\nonumber
\\[-8pt]
\\[-8pt]
\nonumber
&&\qquad \le|E(A_n  \delta_u(G_n))-E(A_n  \delta_u(G)
)|
+ |E(A_n  \delta_u(G))-E(A  \delta_u(G)
)|\\
&&\qquad \le\|A_n\|_{D^{p/(p-1)}_{\a}}\|\delta_u(G_n)-\delta_u(G)\|
_{D^p_{-\a}}
+ \|A_n-A\|_{D^{p/(p-1)}_{\a}}\|\delta_u(G)\|_{D^p_{-\a}},\nonumber
\end{eqnarray}
where we recall Theorem~\ref{theoremwatanabe} for definitions of $p$
and $\a$. We also note that,
since $G_n$ and $G$ are elements of $D^{\infty-}_3$ and they are
nondegenerate, existence of such
$p$ and $\a$ is ensured. Moreover, since $G_n\to G$ in $D^{\infty
-}_3$, it's easy to see that
$\sup_n\|A_n\|_{D^{p/(p-1)}_{\a}} <\infty$ and $\|A_n-A\|
_{D^{p/(p-1)}_{\a}}\to0$. Also,
$\|\delta_u(G_n)-\delta_u(G)\|_{D^p_{-\a}}\to0$, which proves the result.
\end{pf}

\begin{pf*}{Proof of Theorem \protect\ref{theoremGKF} (\textup{Continued})}
We extend $F_n$ from $\real^n$ to $\real^{\infty}$ or, equivalently,
to $X$, by suppressing all the
indices after the first $n$, that is, considering $F_n$ as cylindrical
Wiener functionals.\vspace*{1pt} Then,
by using the invariance property of GMFs, the $\minkmu_j$'s of the extended
$F_n$ remain the same as that of $F_n$ when restricted to $\real^n$.
Together with this, using the
fact that $F_n$ converges to $F$ in $D^{\infty-}_{3}$, we clearly have
%
%e56 #&#
\begin{equation}
\label{eqconvergenceGMF}
\lim_{n\to\infty}\minkmu_j(F_n^{-1}[u,\infty)) = \minkmu
_j(F^{-1}[u,\infty)).
\end{equation}
Therefore, using \eqref{eqGKFn} and \eqref{eqconvergenceGMF}, we
shall have
%
%e57 #&#
\begin{equation}
\label{eqGKFnconvergeI}
\lim_{n\to\infty}E(\lips_0(A_u(f_n;M)))=\sum_{j=0}^{N}c_{j} \lips
_{j}(M) \minkmu_j(F^{-1}[u,\infty)).
\end{equation}

Now, it suffices to prove that the right-hand side of \eqref
{eqexpectationmetan} converges to
the right-hand side of \eqref{eqexpectationmeta}. Clearly,
%
%e58 #&#
\begin{equation}
\label{eqconvergencedensity}
\lim_{n\to\infty}p_{\nabla f_n}(y) = p_{\nabla f}(y),
\end{equation}
which follows from the fact that $\nabla f_n$ converges to $\nabla f$
in a much stronger sense,
as is clear from the assumption $f_n\to f$ in $D^{\infty}_{3+\delta
}$. Next, we need to prove
%
%e59 #&#
\begin{eqnarray}
\label{eqmetaintegrandconvergence}
\qquad &&
\bigl|E\bigl(|\det\nabla^2 f_n(x)|  1_{D_k}(\nabla^2 f_n(x))
1_{[u,\infty)}(f_n(x))|\nabla f_n(x)=0\bigr)
\nonumber
\\[-8pt]
\\[-8pt]
\nonumber
&&\qquad{} -E\bigl(|\det\nabla^2 f(x)|  1_{D_k}(\nabla^2 f(x)) 1_{[u,\infty
)}(f(x))|\nabla f(x)=0\bigr)\bigr| \to0,
\end{eqnarray}
which is similar to the proof of Theorem \ref
{theoremGMFconvergence}. Using precisely the same techniques,
writing $B_n(x) = (|\det\nabla^2 f_n(x)|  1_{D_k}(\nabla^2
f_n(x)) 1_{[u,\infty)}(f_n(x)))$ and defining
$B(x)$ in a similar fashion, we have
\begin{eqnarray*}
&&\bigl|E^{\nabla f_n(x)=0}B_n(x) - E^{\nabla f(x)=0}B(x)\bigr|\\
&&\qquad \le\bigl|E^{\nabla f_n(x)=0}B_n(x) -E^{\nabla f(x)=0}B_n(x)\bigr|\\
&&\qquad\quad{} + \bigl|E^{\nabla
f(x)=0}B_n(x)-E^{\nabla f(x)=0}B(x)\bigr|\\
&&\qquad = |E(B_n(x) \delta_0(\nabla f_n(x))) - E(B_n(x)
\delta_0(\nabla f(x)))|\\
&&\qquad\quad{}  + |E(B_n(x) \delta_0(\nabla f(x))) - E(B(x)
\delta_0(\nabla f(x)))|\\
&&\qquad\le\|B_n(x)\|_{L^{p/(p-1)}}\|\delta_0(\nabla f_n(x)) - \delta
_0(\nabla f(x))\|_{L^p}\\
&&\qquad\quad{}  + \|B_n(x) - B(x)\|_{L^{p/(p-1)}}\|\delta_0(\nabla f(x))\|_{L^p},
\end{eqnarray*}
which, under the assumptions of $f_n(x)\to f(x)$ in $D^{\infty-}_3$
and nondegeneracy of $\nabla f(x)$,
converges to zero as $n\to\infty$,
thus proving that the integrand of \eqref{eqexpectationmetan}
converges to that of
\eqref{eqexpectationmeta} for each $x\in M$.
Then, in order to prove that the integral involved in equation
\eqref{eqexpectationmetan}
converges to the integral in \eqref{eqexpectationmeta}, note that
the random fields $f_n$ and $f$
defined on the manifold $M$ are chosen to be sufficiently smooth so
that we can use an uniform integrability
argument to conclude that the right-hand side of \eqref
{eqexpectationmetan} converges to
the right-hand side of \eqref{eqexpectationmeta}. Therefore, we
shall have
\begin{eqnarray*}
E(\lips_0(A_u(f;M))) & = & \lim_{n\to\infty} E(\lips
_0(A_u(f_n;M)))\\
& = & \lim_{n\to\infty}\sum_{j=0}^m c_j \lips_j(M) \minkmu
_j(F_n^{-1}[u,\infty))\\
& = & \sum_{j=0}^m c_j \lips_j(M) \minkmu_j(F^{-1}[u,\infty)),
\end{eqnarray*}
where in going from the first line to the second, we have used the
finite dimensional results set forth in
\cite{RFG}, and in going from the second to the third line, we have
used Theorem~\ref{theoremGMFconvergence}.
\end{pf*}

%%-------------------------------------------------------------------------------------------------------------------------------------------------------------------------
%%-------------------------------------------------------------------------------------------------------------------------------------------------------------------------

% imsref loaded by akundreckaite, 2012-05-23 14:44:36

%suskaldyti doi

\printaddresses

\end{document}